\documentclass[aop,MSNbibl,nameyear,dvips]{arximspdf}

\doi{10.1214/12-AOP778} 
\volume{41}
\issue{5}
\pubyear{2013}
\firstpage{3617}
\lastpage{3657}

\makeatletter

\newcommand{\iint}{\int\!\!\int}

\newproclaim{Remarks}{Remarks}

\newtheorem{theorem}{Theorem}
\newtheorem{lemma}[theorem]{Lemma}

\newcommand{\N}{{\mathbb N}}

\makeatother

\begin{document}
\begin{frontmatter}

\title{Conditioning super-Brownian motion on its boundary statistics, and
fragmentation}
\runtitle{Conditioning super-Brownian motion}

\begin{aug}
\author[A]{\fnms{Thomas S.} \snm{Salisbury}\thanksref{t1}\ead[label=e1]{salt@yorku.ca}}
\and
\author[B]{\fnms{A. Deniz} \snm{Sezer}\corref{}\thanksref{t2}\ead[label=e2]{adsezer@ucalgary.ca}}
\runauthor{T. S. Salisbury and A. D. Sezer}
\affiliation{York University and University of Calgary}
\address[A]{Department of Mathematics\\
\quad and Statistics\\
York University\\
4700 Keele St.\\
Toronto, ON M3J 1P3 \\
Canada\\
\printead{e1}} 
\address[B]{Department of Mathematics\\
\quad and Statistics\\
University of Calgary\\
2500 University Drive NW\\
Calgary, AB T2N 1N4 \\
Canada\\
\printead{e2}}
\end{aug}

\thankstext{t1}{Supported in part by NSERC.}

\thankstext{t2}{Supported in part by NSERC.
Much of this research took place during D. Sezer's post-doctoral visit
at York University, 2005--2008.}

\received{\smonth{10} \syear{2011}}
\revised{\smonth{5} \syear{2012}}

%
\begin{abstract}
We condition super-Brownian motion on ``boundary statistics'' of the
exit measure $X_D$ from a bounded domain $D$. These are random
variables defined on an auxiliary probability space generated by
sampling from the exit measure $X_D$. Two particular examples are:
conditioning on a Poisson random measure with intensity $\beta X_D$ and
conditioning on $X_D$ itself. We find the conditional laws as
$h$-transforms of the original SBM law using Dynkin's formulation of
$X$-harmonic functions. We give explicit expression for the (extended)
$X$-harmonic functions considered. We also obtain explicit
constructions of these conditional laws in terms of branching particle
systems. For example, we give a fragmentation system description of
the law of SBM conditioned on $X_D=\nu$, in terms of a particle
system, called the backbone. Each particle in the backbone is labeled
by a measure $\tilde{\nu}$, representing its descendants' total
contribution to the exit measure. The particle's spatial motion is an
$h$-transform of Brownian motion, where $h$ depends on $\tilde{\nu}$.
At the particle's death two new particles are born, and $\tilde{\nu}$
is passed to the newborns by fragmentation.
\end{abstract}

%

\end{frontmatter}

\section{Introduction}\label{introsection}

Studying conditioned Markov processes is a kind of inverse
problem---given information about how the process ends up, one tries to
infer how it got there, at least in terms of probabilities. In the
context of Brownian motion and finite-dimensional Markov processes, one
can make very explicit calculations, starting with the work of
\citet{Doob}. Attempts to make similar calculations for super
Brownian motion are more recent. These studies typically aim to recover
the conditional law of a superprocess as the law of a distinct
probabilistic object. Several authors have succeeded in coming up with
such descriptions for certain conditionings, and produced models with
remarkably rich structure. The first of these models was the immortal
particle system of \citet{EP90} and \citet{Eva93}, for
super-Brownian motion in $\mathbb{R}^{d}$ conditioned on survival. In
this model, an immortal particle moves according to a Brownian motion,
and throws off mass at a uniform rate, and then this mass evolves in
the space as an unconditioned super-Brownian motion. Following Evans
and Perkins, \citet{SV} considered a super Brownian motion $X$ in
a domain $D$, with conditioning based on the \textit{exit measure} $X_D$
from $D$. More specifically, they conditioned $X$ on the event that the
support of $X_{D}$ contains certain points $z_1,\ldots,z_k$, and
recovered the resulting conditional law in terms of a branching
backbone system. The branching backbone is a random tree with $k$
leaves reaching the points $z_1,\ldots, z_k$. Similar to Evans and
Perkins's model, there is mass uniformly created along the branching
backbone which follows the law of an unconditioned super-process
independent of the points $z_1,\ldots,z_k$. Giving such an explicit
characterization of a conditioned process is an interesting problem
from a probabilistic modeling point of view. For example, in population
dynamics, one can view it as an analogue of a host of biological
problems in which one has information about the state of the population
at certain times or locations, and one wishes to infer the genealogical
structure of the populations of the ancestors (e.g., the ``out of
Africa'' problem of human origins). For explicit representations of
other related conditioned processes, see \citet{RR89},
Overbeck (\citeyear{O93,O94}), \citet{E93} or \citet{EO03}.

It turns out that there is more to the conditioning problem than described
above. A conditioned process represents a special case of a Girsanov
transformation, or a martingale change of measure. For concreteness,
let us consider the following example: let $\xi$ be Brownian motion in
a domain $E$. Let $\tau_E$ be the exit time from $E$. We compute the
conditional law, $\Pi^{z}_x$ of $\xi$ given $\xi_{\tau_E}=z$, by a
martingale change of measure from $\Pi_x$, the law of $\xi$. This
martingale change of measure is given in terms of a certain harmonic
function $h^{z}(\cdot)$.
More precisely, for any domain $D$ compactly contained in $E$, and any
$Y$ measurable with respect to $\mathcal{F}_{\tau_{D}}$, we have $\Pi
^{z}_x(Y)=\Pi_x(Y h^{z}(\xi_{\tau_{D}})/h^z(x))$.
A Girsanov transformation defined in terms of a harmonic function is
called an $h$-transform, and typically conditional laws of Markov
processes are formulated as $h$-transforms of their original laws. This
classical relationship between harmonic functions and conditioning a
Markov process leads to an elegant probabilistic formulation of the
Martin boundary theory for elliptic differential operators.

In the context of super-processes, the analogue of harmonic functions
are $X$-harmonic functions. Following
the definition of \citet{Dyn02},
let us consider a super-Brownian motion $X=(X_{D}, P_\mu)$, a family of
random measures (\textit{exit measures}) and their associated probability
laws where $D$ is an open subset of a given domain $E$ in $\mathbb
{R}^{d}$, and $\mu$ is a finite measure on $E$. Write $D\Subset E$ if
$D$ is open and its closure is a compact subset of $E$. A nonnegative
function $H$ is \textit{$X$-harmonic} if for any $D\Subset E$ and any
finite measure $\mu$ with support in $D$,
%
%
\begin{equation}\label{eqmvp}
P_{\mu}\bigl(H(X_D)\bigr)=H(\mu).
\end{equation}
Note that this property resembles the mean value property of a
harmonic function, hence the name $X$-harmonic.
Moreover, the $X$-harmonic functions are related to conditioning
super-Brownian motion in the same way as harmonic functions are related
to conditioning Brownian motion; they give us the explicit Girsanov
transformation to switch from the unconditioned probability law to the
conditioned probability law. An
$H$-transform $P_{\mu}^{H}$ is obtained from $P_\mu$ by setting
\[
P_{\mu}^{H}(Y)=\frac{1}{H(\mu)}P_{\mu}
\bigl(H(X_{D})Y\bigr)
\]
for $Y$ nonnegative and $\mathcal{F}_{\subset D}$-measurable, where
$\mathcal{F}_{\subset D}$ is the $\sigma$-algebra generated by $X_{D'},
D'\subset D$.
In his book, \citet{Dyn02}
suggests a new direction for investigating the solutions of the p.d.e.
$\frac12\Delta=2u^{2}$, namely to explore $X$-harmonic functions
(thinking them as the analogue of harmonic functions) and ultimately,
to build a Martin boundary theory for this nonlinear p.d.e. In this
case, the Martin boundary is defined as the set of extreme elements of
the convex set of all $X$-harmonic functions. Dynkin points out that a
concrete understanding of extreme $X$-harmonic functions might yield
further insights about the solutions of the p.d.e. $\frac12\Delta u=2u^{2}$.
Since then major progress has been made on the study of the solutions
of this p.d.e. using other approaches. For example, \citet{M04}
classified the solutions as $\sigma$-moderate in the case of a smooth
domain (solving a conjecture of Dynkin and Kuznetsov). However, the
relationship between \mbox{$X$-harmonic} functions and solutions remains
largely unexplored. Dynkin, in a series of papers, has taken concrete
steps to better formulate and understand extreme $X$-harmonic
functions. \citet{Dyn06} obtains the extreme $X$-harmonic
functions by a
limiting procedure from the Radon--Nykodym densities $H^{\nu}_{D}(\mu
)=\frac{dP_{\mu,X_D}}{dP_{c,X_D}}(\nu)$ of $P_{\mu,X_D}(d\nu
)=P_\mu
(X_D\in d\nu)$ with respect to $P_{c,X_D}(d\nu)=P_{c}(X_D\in d\nu)$.
\citet{DynkinIJM06} derives a formula for $H^{\nu}_{D}(\mu)$ using
diagram description of moments. The functions $H^{\nu}_{D}$ will be
central to our analysis as well; we will run into them while studying
conditionings of SBM, and use them to derive results about the
structure of conditioned SBM.

Our goal in this paper is to explore various ways of conditioning a
super-Brownian motion. We are motivated by the rich structure of the
underlying probabilistic objects as well as its potential connection to
Dynkin's research program on Martin boundary theory of SBM. Here is a
summary of our contributions: we develop a
way of conditioning a super-Brownian motion, which we call
``conditioning on boundary statistics.'' The random variables which we
condition on are defined on an auxiliary probability space, and
generated by sampling from the exit measure $X_D$. We find
representations of these conditionings as $H$-transforms of
unconditioned SBM. In general, we identify the functions $H$ as
``extended'' $X$-harmonic. The term ``extended'' is used because, even
though we show that these functions satisfy the mean value property
(\ref{eqmvp}), in general, we do not know whether they are finite for
all $\mu$. An example of a boundary statistic is
a Poisson random measure with intensity $\beta X_D$, where $\beta>0$.
It turns out that for this kind of conditioning the resulting
$H$-transform is through an $X$-harmonic family of functions studied
earlier in the literature. Another important and more complex example
we studied is conditioning SBM on its exit measure $X_D$, which is what
most of our paper is devoted to.
We find that the corresponding $H$-transform uses the family $H^{\nu
}_{D}(\mu)=\frac{dP_{\mu,X_D}}{dP_{c,X_D}}(\nu)$, densities first
introduced in \citet{Dyn06}. This paper shows that one can choose a
version of this family such that for each $\nu$, $H^{\nu}(\mu)$ will
satisfy the mean value property (\ref{eqmvp}) for each $\mu$. Although
this family is claimed to be $X$-harmonic in \citet{Dyn06}, we are not
aware of any results that actually show that $H^{\nu}$ will be a finite
function of $\mu$. In this paper we shall classify $H^{\nu}_{D}$ as an
extended $X$-harmonic function and leave the question of finiteness to
be resolved in a different paper, as this will require us to develop
analytical bounds on the densities of moment measures of SBM. In this
paper we will not go beyond describing the probabilistic structure of
conditioned SBM on its exit measure; as such we are
not going to lose much generality by stating and proving our results
without assuming that $H^{\nu}$ is finite. Indeed, for our purposes it
will suffice that $H^{X
_D}_{D}(\mu)$ is finite $P_{\mu}$ almost surely, which is true. The
family $H^{\nu}_{D}(\mu)$ is of special interest because it
can be considered as the analogue of the Poisson kernel of $D$. Also,
\citet{Dyn06}
showed that if $H$ is an extreme $X$-harmonic function in $E$, then for
every $\mu$, and for every sequence $D_k$ exhausting $E$
\[
H(\mu)=\lim_{k\rightarrow\infty}H_{D_k}^{X_{D_k}}(\mu)
\]
$P^{H}_\mu$ almost surely.

The heart of the paper is Theorem~\ref{theolinear}, which gives a new
formula for $H_{D}^{\nu}$. From this we deduce an infinite particle
fragmentation system description of $P_{\mu}^{\nu}=P_{\mu}^{H^\nu}$,
the conditional law of $X$ in $D$ given $X_{D}=\nu$ (Theorem \ref
{theofrag}). This is carried out in terms of a particle system, called
the backbone in \citet{SV}, along which a mass is created uniformly.
In the backbone, each particle is assigned a measure $\tilde{\nu}$ at its
birth. The spatial motion of the particle is an $h$-transform
of\vadjust{\goodbreak}
Brownian motion, where $h$ is a potential that depends on $\tilde{\nu
}$. The measure
$\tilde{\nu}$ represents the particle's
contribution to the exit measure. At the particle's death two new
particles are born, and $\tilde{\nu}$ is passed to the newborns by
fragmentation into two bits. Here, we use the techniques of \citet{SV}
applied to a more general setting. This description connects the theory
of conditioned super-processes to the growing literature on infinite
fragmentation and coalescent processes; see, for example,
\citet{Bertoin} for a comprehensive exposition.

\section{Preliminaries}\label{prelimsection}
\subsection{Super-Brownian motion}

We will follow Dynkin's definition of super-Brownian motion (SBM). Let
$E$ be a domain of $\mathbb{R}^{d}$, and let $\mathcal{M}_E$ be the
positive finite measures on $E$.
A super-Brownian motion, $(X_{D}, P_\mu)$, is a family of random
measures (\textit{exit measures}) and their associated probability laws
where $D$ is an open subset of a given domain $E$ in $\mathbb{R}^{d}$,
and $\mu$ is a finite measure on $E$ with the following properties:
\begin{longlist}[(a)]
\item[(a)] Exit property: $P_{\mu}(X_{D}(D)=0)=1$
for every $\mu$, and if $\mu(D)=0$, then $P_{\mu}(X_{D}=\mu)=1$.
\item[(b)] Markov property: if
$Y\ge0$ is measurable with respect to the $\sigma$-algebra $\mathcal
{F}_{\subset D}$ generated by $X_{D'},D'\subset D$ and $Z\geq0$ is
measurable with respect to the $\sigma$-algebra $\mathcal{F}_{\supset
D}$ generated by $X_{D''}$, $D''\supset D$, then
\[
P_\mu(YZ)=P_\mu(YP_{X_{D}}Z).
\]
\item[(c)] Branching property: for any nonnegative Borel $f$,
\[
P_\mu
\bigl(e^{-\langle X_{D},f\rangle}\bigr)= e^{-\langle\mu, V_Df\rangle}\qquad \mbox{where }
V_Df(y)=-\log P_{y}\bigl(e^{-\langle X_D,f\rangle}\bigr)
\]
and $P_y=P_{\delta_y}$.
\item[(d)] Integral equation for the log-Laplace functional: $V_{D}f$ solves
the integral equation
\[
u+G_D\bigl(2u^2\bigr)=K_D f,
\]
where $G_D$ and $K_D$ are, respectively, Green and Poisson operators for
Brownian motion in $D$. In other words, if $\xi_t$ is a Brownian motion
starting from $x$, under a probability measure $\Pi_{x}$,
then $K_{D}f(x)=\Pi_xf(\xi_{\tau_D})$, where $\tau_D$ is the exit time
from $D$. Likewise, $G_{D}f(x)=\Pi_x(\int_{0}^{\tau_D}f(\xi_t)\,dt)$.
\end{longlist}

Under certain regularity conditions on $D$ and $f$ [see, e.g.,
\citet{Dyn02}], the integral equation in (d) is equivalent to the boundary
value problem
\begin{eqnarray*}
\tfrac12\Delta u&=&2u^2,
\\
u(x)&=&f(x),\qquad x\in\partial D.
\end{eqnarray*}

$X_D$ represents the exit measure from $D$, and the first property
simply means that $X_D$ is concentrated on $D^c$, so that exiting is
instantaneous if we start outside $D$.\vadjust{\goodbreak} The third property means that
distinct clumps of initial mass evolve independently. It follows from
continuity of Brownian motion that $X^D$ is supported on $\partial D$,
if the initial measure is supported on $D$; see Property 2.2.A of \citet{Dyn02}.
The fourth property restricts attention to finite variance branching,
and normalizes the branching rate. The normalizing factor 2 in front of
$u^2$ is chosen to be consistent with \citet{LG} and \citet{SV}.

\subsection{Infinite divisibility and Poisson representation}
It is well known that $X_{D}$ has an infinitely divisible distribution
for each $D$. This property leads to the construction of a new measure,
$\mathbb{N}_{x}$, called the super-Brownian excursion law starting from
$x$. Under $\mathbb{N}_{x}$, $X$ evolves as a super-Brownian motion,
but $\mathbb{N}_x$ will be $\sigma$-finite, not a probability. Thus, it
is technically more complicated than $P_{\mu}$. But $\mathbb{N}_{x}$ is
actually a more basic object heuristically, under which the genealogies
are simpler, because
all the mass starts from a single particle located initially at $x$.

In fact $P_{\mu}$ can be built up from a Poisson random measure with
intensity $\theta(d\chi)=\int\mathbb{N}_{x}(d\chi)\mu(dx)$. More
precisely, let
\[
\Pi(d\chi)=\sum_{i}\delta_{\chi^{i}}
\]
be such a Poisson random measure, where the $\chi^{i}$ are random
measure valued paths. Then $X=\sum\chi^i=\int\chi\Pi(d\chi)$
is a super-Brownian motion with initial state~$\mu$. In terms of
$X_{D}$, this yields the following formula; see Theorem 5.3.4 of
\citet{Dyn04a}: let $F$ be a nonnegative measurable function
defined on
$\mathcal{M}_{E}$. Then
%
%
\begin{eqnarray}
\label{eqPoissonrep}
&&
P_{\mu}\bigl(F(X_{D})\bigr)\nonumber\\
&&\qquad=
e^{-\mathcal{R}_\mu(\mathcal{M}_E)}F(0)
\\
&&\qquad\quad{}+\sum_{n=1}^{\infty}\frac{1}{n!}e^{-\mathcal{R}_\mu(\mathcal
{M}_E)}
\int\mathcal{R}_\mu(d\nu_1)\cdots\mathcal{R}_\mu(d
\nu_n)F(\nu_1+\cdots+\nu_n),\nonumber
\end{eqnarray}
where $\mathcal{R}_\mu$ is the canonical measure of $X_{D}$ with
respect to $P_\mu$ and can be derived from $\N_x$ by $\mathcal
{R}_\mu
(A)=\langle\mu, \N_{(\cdot)}(X_D\in A,X_D\neq0)\rangle$.
In other words, we obtain $X_D$ as a superposition of a Poisson number
of more ``basic'' exit measures, each descended from a single initial
individual. These ``basic'' exit measures arise as the atoms of a
Poisson random measure whose characteristic measure is $\mathcal
{R}_\mu$.

Note that a special case of the above representation gives us
$V_{D}f(x)=\N_x(1-e^{-\langle X_D, f\rangle})$. Note further that we
will, in the future, write this as
\[
P_{\mu}\bigl(F(X_{D})\bigr)= \sum
_{n=0}^{\infty}\frac{1}{n!}e^{-\mathcal{R}_\mu(\mathcal
{M}_E)}\int
\mathcal{R}_\mu(d\nu_1)\cdots\mathcal{R}_\mu(d
\nu_n)F(\nu_1+\cdots+\nu_n)
\]
by taking the convention that the $n=0$ term is the first expression in
(\ref{eqPoissonrep}).\vadjust{\goodbreak}

The measure $\N_x$ was first considered by
\citet{LG}
in his random snake formulation of super-Brownian motion. As we follow
Dynkin's framework to study our problem, we refer the reader to
\citet
{Dyn04a} for a systematic account of the theory of the measures $\N_x$
and their applications. Note that the latter has a general branching
function $\psi$, but for us, this is taken to be $\psi(u)=2u^2$.

\subsection{Moment measures of super Brownian motion} Among the key
tools in our analysis are the recursive moment formulas of SBM. The
moment measures of SBM are the following measures:

Let $\phi,f_1,\ldots,f_n$ be positive Borel functions and write
$f=(f_1,\ldots,f_n)$.
For
$C\subset\{1,\ldots,n\}$, let
%
%
\begin{eqnarray}
\label{defnmoment}
n_C(\phi,f,x)&=&\N_{x}e^{-\langle X_D,\phi\rangle}
\Pi_{i\in
C}\langle X_D,f_i\rangle,
\\
\label{defpmoment1}
p_C(\phi,f,\mu)&=&P_{\mu}e^{-\langle X_D,\phi\rangle}
\Pi_{i\in
C}\langle X_D,f_i\rangle.
\end{eqnarray}
Let $K_{D}^{l}$ and $G_{D}^{l}$ be the Poisson and Green operator
for the operator
\[
\mathcal{L}^{l}=\tfrac{1}{2}\triangle-l,
\]
where $l(x)=4V_D \phi(x)$. In other words, let $\xi_t$ be a diffusion
starting from $x$, with generator $\mathcal{L}^{l}$ under a probability
measure $\Pi_{x}^{l}$.
Then $K_{D}^{l}f(x)=\Pi_x^{l}f(\xi_{\tau_D})$ and $G_{D}^{l}f(x)=\Pi
_x^l(\int_{0}^{\tau_D}f(\xi_t)\,dt)$, where $\tau_D$ is the exit time
from~$D$.

For $C=\{i\}$ we have the Palm formula
%
%
\begin{equation}\label{eqpalm1}
n_C(\phi,f,x)=K^{l}_{D}f_i(x),
\end{equation}
and for general $C$ we have the following recursive formulas; see, for
example, Theorem 5.1.1 of \citet{Dyn04a}, or Lemma 2.6 of
\citet{SV}:
%
%
\begin{eqnarray}
\label{eqrecur1}
n_C(\phi,f,\cdot)&=&\frac{1}{2}\sum
_{A \subset C, A\neq\varnothing,C}G_{D}^{l}\bigl(4n_A(
\phi,f,\cdot)n_{C\setminus
A}(\phi,f,\cdot)\bigr),
\\
\label{eqrecur2}
p_C(\phi,f,\mu)&=&e^{-\langle
\mu,V_D(f_0)\rangle}\sum
_{\pi(C)}\bigl\langle\mu,n_{C_1}(\phi,f,\cdot)\bigr
\rangle\cdots\bigl\langle\mu,n_{C_r}(\phi,f,\cdot)\bigr\rangle.
\end{eqnarray}
Here $\pi(C)$ is the set of partitions of $C$. These formulas will
allow us to construct a variety of extended $X$-harmonic functions of
polynomial type.

We will also need extensions of these formulas for
\[
\N_{x}e^{-(\langle X_{D_1},\phi_1\rangle+\cdots+\langle X_{D_k},\phi_k
\rangle)}\Pi_{i\in C}\langle
X_D,f_i\rangle
\]
and
\[
P_{\mu}e^{-(\langle X_{D_1},\phi_1\rangle+\cdots+\langle
X_{D_k},\phi_k
\rangle)}\Pi_{i\in C}\langle X_D,f_i
\rangle,
\]
where $D_i\subset D_k=D$. Formulas (\ref{eqpalm1}), (\ref{eqrecur1})
and (\ref{eqrecur2}) give us these quantities when $k=1$ and $x\in D$.
For $k\geq2$ we find them recursively as follows. Let us put $I= \{
D_1,D_2,\ldots,D_k\}$, $\phi_{I}=(\phi_1,\ldots,\phi_k)$,
$u^I(x)=\mathbb{N}_{x}(1-\exp-(\langle X_{D_1},\phi_1\rangle+\cdots
+\langle X_{D_k},\phi_k \rangle))$, $l^I=4u^I$, $D_I=D_1\cap\cdots
\cap
D_k$ and $I_j=I-\{D_j\}$. We define an operator $n_C^I(\phi_I,f,x)$ as
follows. For $\operatorname{card}(I)=1$,
%
%
\begin{equation}
\label{defnmoment1} n_C^I(\phi,f,x)=\cases{
\N_{x}e^{-(\langle X_{D},\phi\rangle)}\Pi_{i\in
C}\langle X_D,f_i
\rangle, &\quad$x\in D$,
\cr
f_i(x), &\quad$x\notin D, C=\{i\}$,
\cr
0, &\quad$x
\notin D, \operatorname{card}(C)>1$,}
\end{equation}
and for $\operatorname{card}(I)>1$,
%
%
\begin{eqnarray}
\label{defnmoment2}
&&n_C^I(\phi_I,f,x)\nonumber\\[-8pt]\\[-8pt]
&&\qquad=
\cases{\N_{x}e^{-(\langle X_{D_1},\phi_1\rangle+\cdots
+\langle
X_{D_k},\phi_k \rangle)}\Pi_{i\in C}\langle
X_D,f_i\rangle, &\quad$x\in D_I$,
\vspace*{3pt}\cr
n_C^{I_j}(\phi_{I_j},f,x), &\quad$x\notin
D_j, j\in I$.}\nonumber
\end{eqnarray}
We let
%
%
\begin{equation}
\label{defpmoment2}\quad
p_C^I(\phi_1,\phi_{2},\ldots,
\phi_k, f,\mu)=P_{\mu}e^{-(\langle
X_{D_1},\phi_1\rangle+\cdots+\langle X_{D_k},\phi_k \rangle)}\Pi_{i\in
C}
\langle X_{D_k},f_i\rangle.
\end{equation}
Fix an $f$ and write $n_C^I(x)$ for $n_C^I(\phi_I, f,x)$ and
$p_C^I(\mu
)$ for $p_C^I(\phi_I,f,\mu)$. The following formulas tell us that the
values of $n_C^I(x)$ on $D_I$ can be recovered from its values at the
boundary of $D_I$ and the values of the functions $n_{A}^{I}(x)$ on
$D_I$, with $A\subset C$ and $A\neq\varnothing, C$, using the Poisson
and Green operators of $\mathcal{L}^{l^I}$ for the domain~$D_I$.
For $C=\{i\}$ and $\operatorname{card}(I)>1$,
%
%
\begin{equation}
\label{eqexpalm1} n_C^I(x)=K^{l^I}_{D_I}
\bigl(n_C^{I}\bigr) (x),\qquad x\in D_I,
\end{equation}
and for $\operatorname{card}(C)>1$ and $\operatorname{card}(I)>1$,
%
%
\begin{equation}\label{eqexrecur1}\qquad
n_C^I(x)=\frac{1}{2}\sum
_{A \subset C, A\neq\varnothing,C}G_{D_I}^{l^I}\bigl(4n_A^In_{C\setminus
A}^I
\bigr) (x)+K^{l^I}_{D_I}\bigl(n_C^I
\bigr) (x),\qquad x\in D_I.
\end{equation}
Then for $\mu$ compactly supported in $D_I$
%
%
\begin{equation}\label{eqexrecur2}
p_C^I(\mu)=e^{-\langle
\mu,l^I\rangle}\sum
_{\pi(C)}\bigl\langle\mu,n^I_{C_1}\bigr
\rangle\cdots\bigl\langle\mu,n^{I}_{C_r}\bigr\rangle.
\end{equation}

Since at the boundary of $D_I$, $n_C^I(x)$ is recursively defined in
terms of $n_C^{I_j}(x)$, we have now a complete recursive algorithm to
compute $n_C^I(x)$, starting with $I=\{D\}$, and $C=\{i\}$, recursively
first building $n_C^I(x)$ for all $C$ keeping $I$ the same, and then
increasing the cardinality of $I$ by 1 and repeating the same procedure
until the desired cardinality of $I$ is achieved.

We omit the proof of these formulas and refer the reader to the proof
of Theorem 5.1.1 of \citet{Dyn04a}. The reader will realize that the
argument in that proof works also for the functions $u_I^{\lambda}(x)$
which\vadjust{\goodbreak} are defined as follows. Set $X_I=(X_{D_1},\ldots,X_{D_k})$, so
$\langle X_I,\phi_I\rangle=\sum_j\langle X_{D_j},\phi_j\rangle$. For
for $I=\{D\}$, set
\[
u^\lambda_I(x)=\cases{\mathbb{N}_x\bigl(1-
\exp{-\bigl(\langle X_{I},\phi_I\rangle+
\lambda_1\langle X_{D_k},f_1 \rangle+\cdots+
\lambda_n\langle X_{D_k},f_n\rangle\bigr)}
\bigr),
\cr
\hphantom{\phi_k(x)+\lambda_1f_1(x)+\cdots+\lambda_nf_n(x),\,}
\qquad x\in D,
\cr
\phi_k(x)+\lambda_1f_1(x)+
\cdots+\lambda_nf_n(x), \qquad x\notin D,}
\]
and recursively for $\operatorname{card}(I)>1$ as
\[
u^\lambda_I(x) = \cases{\mathbb{N}_x\bigl(1-
\exp{-\bigl(\langle X_{I},\phi_I\rangle+
\lambda_1\langle X_{D_k},f_1 \rangle+\cdots+
\lambda_n\langle X_{D_k},f_n\rangle\bigr)}
\bigr),
\cr
\hphantom{u^{\lambda}_{I_j}(x)+\phi_j(x),\,}
\qquad x\in D_I,
\cr
u^{\lambda}_{I_j}(x)+\phi_j(x),
\qquad x\notin D_j,j\in I,j\neq k.}
\]
We note that $u_{I}^{\lambda}$ satisfies for $x\in D_I$,
\[
u_{I}^{\lambda}+G_{D_I}^{l^I} \bigl(2
\bigl(u_{I}^{\lambda
}\bigr)^2\bigr)=K_{D_I}^{l^I}
\bigl(u_{I}^{\lambda}\bigr),
\]
following formula 2.11 of \citet{Dyn02}, Chapter 3, which in turn yields
the formulas for $n_C^I$ by differentiating $u_{I}^{\lambda}$ with
respect to $\lambda$.

\subsection{Absolute continuity}
The moment formulas together with the Markov property and Poisson
representation yield an important theorem taken in this form from
Theorem 5.3.2 of \citet{Dyn04a}. See also Proposition 2.18 of
\citet{M04}. Let $\mathcal{M}^{c}_{D}$ be the space of finite measures
compactly supported in $D$.
%
%
\begin{theorem}\label{theoabscon} Suppose $A\in\mathcal{F}_{\supset
D}$. Then either $P_\mu(A)=0$ for all $\mu\in\mathcal{M}^{c}_{D}$ or
$P_\mu(A)>0$ for all $\mu\in\mathcal{M}^{c}_{D}$.
\end{theorem}

\section{Extended $X$-harmonic functions and conditioning}
\label{Xharmonicsection}
\label{secboundstat}

In the rest of the paper we fix $D\Subset E$. Following
\citet{Dyn06},
a nonnegative function $H\dvtx\mathcal{M}^c_D\to[0,\infty)$ is called
\textit{$X$-harmonic} in $D$, if for any $D'\Subset D$ and any finite
measure $\mu\in\mathcal{M}^c_{D'}$,
%
%
\begin{equation}
\label{eqmvp2} P_{\mu}\bigl(H(X_{D'})\bigr)=H(\mu).
\end{equation}
%
We will call a nonnegative $H$ extended $X$-harmonic if it satisfies
(\ref{eqmvp2}) but is not necessarily everywhere finite.

We are going to touch upon three different kinds of extended
$X$-harmonic functions, which are derived from conditioning SBM on
its various boundary statistics. These boundary statistics are:
\begin{longlist}[(a)]
\item[(a)] a Poisson random measure with characteristic measure
$\beta X_D$;
\item[(b)] a random variable $Z$ drawn from the probability distribution
$\frac{X_{D}}{{\langle X_D,1\rangle}}$ if \mbox{$X_D\neq0$}, and set equal
to some given
$\Delta
\notin\partial D$ if $X_D=0$;
\item[(c)]${L}(X_D)$, where ${L}$ is a linear map from $\mathcal
{M}_{\partial D}$ to a vector space $V$ [e.g., ${L}(\mu)=\mu$ or
${L}(\mu)=\langle\mu,1\rangle$ means we condition on $X_{D}$ or on its
total mass].
\end{longlist}
Let $S$ be any one of the above statistics. Let $\Sigma$ be the state
space of $S$. We will assume that $\Sigma$ is endowed\vadjust{\goodbreak} with a countably
generated $\sigma$-algebra $\mathcal{S}$ such that $(\Sigma,\mathcal
{S})$ is a measurable Luzin space; see \citet{Dyn06} for a definition.
For example, when $S=X_D$, $\Sigma$ is $\mathcal{M}_{\partial D}$, the
space of finite measures on $\partial D$ and $\mathcal{S}$ is the
$\sigma$-algebra in $\mathcal{M}_{\partial D}$ generated by the
functions $f(\nu)=\nu(B)$, where $B$ is a Borel subset of $\partial D$.
Given $X_D=\nu$, we let $P_S^\nu$ denote the conditional distribution
of $S$. For example, $P_S^\nu(f)$ equals $\langle
\nu,f\rangle/\langle\nu,1\rangle$ in the second case (provided
$\nu\neq0$), and $f({L}(\nu))$ in the third.

$P_\mu$ denotes a probability measure in which $X$ is an SBM started
from $\mu\in\mathcal{M}_{D}^{c}$, and in which $S$ is then drawn (if
necessary) by further sampling. In other words, $P_\mu$ is a
probability defined on the $\sigma$-field $\mathcal{G}=\mathcal
{F}_{\subset D}\lor\sigma\{S\}$. When $\mu=\delta_x$ we set $P_\mu
=P_x$. By construction, $P_\mu(f(S)\mid\mathcal{F}_{\subset
D})=P_S^{X_D}(f)$. In other words, for any $\mathcal{F}_{\subset
D}$-measurable $Y$ we have that
\[
P_\mu\bigl(f(S)Y\bigr)= P_\mu\bigl(P_S^{X_D}(f)Y
\bigr).
\]
Likewise, we let $P_{\mu,S}$ and $P_{x,S}$ denote the marginal
distribution of $S$ under $P_\mu$ and $P_x$, so $P_{\mu,S}(f)=P_\mu(f(S))$.

Let $\mathcal{F}_{\subset D-}=\sigma\{X_{D'},D'\Subset D\}$.
What we want is the conditional law of $\{X_{D'},D'\Subset D\}$ given
$S=s$, which should
therefore be a transition kernel $P^s_\mu$ from $\Sigma$, the state
space of $S$, to the $\mathcal{F}_{\subset D-}$ measurable functions.
More precisely, we will have $P^{S}_{\mu}(Y)=P_{\mu}(Y|S)$, $P_{\mu}$
a.s. for all $\mathcal{F}_{\subset D-}$ measurable $Y$. The following
theorem tells us how we can construct this transition kernel [part
(d)]. The first three statements [(a), (b) and (c)] of this theorem are
equivalent to Theorem 1.1 of \citet{Dyn06} in the case $S=X_D$.
For a
general $S$, we follow Dynkin's proof, with some modifications.

Let us fix a point $x\in D$.
%
%
\begin{theorem} \label{theocondprob} There exists a family of
nonnegative functions $\{H^{s}_{x}\dvtx\mathcal{M}_{D}^{c}\mapsto R_+,
s\in\Sigma\}$ with the following properties:
\begin{longlist}[(a)]
\item[(a)]$H^{(\cdot)}_{x}(\cdot)\dvtx(s,\mu)\mapsto H^{s}_{x}(\mu)$ is
measurable and strictly positive;
\item[(b)] For all $\mu$, $H^{(\cdot)}_{x}(\mu)\dvtx s \mapsto H^{s}_{x}(\mu
)$ is
a version of $\frac{dP_{\mu,S}}{dP_{x,S}}$;
\item[(c)] For all $s$, $H^{s}_{x}(\cdot)$ is extended $X$-harmonic in $D$;
\item[(d)] Define a probability $P^{s}_{\mu}$ on $\mathcal{F}_{\subset D-}$
by setting
%
%
\begin{equation}\label{eqhtrans}
P^{s}_{\mu}(Y)=\frac{1}{H_{x}^{s}(\mu)}P_{\mu}
\bigl(YH_{x}^{s}(X_{D'})\bigr)
\end{equation}
for all $D' \Subset D$ containing the support of $\mu$, and $\mathcal
{F}_{\subset D'}$-measurable $Y$ [whenever $H_{x}^{s}(\mu)<\infty$,
and otherwise setting $P^{s}_{\mu}$ to an arbitrary probability
measure]. Then $P^{S}_{\mu}$ is a version of the conditional law of $X$
given $S$ with respect to $P_\mu$ for all $\mu$.
\end{longlist}
Any two families satisfying the above properties will coincide for
$P_{x,S}$-a.e. $s\in\Sigma$.
\end{theorem}
\begin{pf} Existence of a family $\{\bar{H}_{x}^{s}, s\in\Sigma\}$
with the first two properties follows from
Theorem A.1 of \citet{Dyn06} and the absolute continuity of the family
$\{P_{\mu,S},\mu\in\mathcal{M}^{c}_{D}\}$ with respect to $P_{x,S}$.
Let $O$ be a subdomain compactly contained in $D$. Then
%
%
\begin{equation}\label{eqDynkin}
P_{\mu}\bar{H}^{s}_{x}(X_O)=
\bar{H}_{x}^{s}(\mu)
\end{equation}
for $P_{x,S}$-a.e. $s$, $\forall\mu\in\mathcal{M}^{c}_{O}$. \citet
{Dyn06} proves this
when $S=X_{D}$, and the proofs for the other cases are almost
identical to his.
Next, we want to construct an extended $X$-harmonic function $H^{s}$
for all $s\in\Sigma$.
To do this, we choose a
countable base $O_n$ (w.l.o.g. closed under finite unions), and
probability measures
$\mu_n\in\mathcal{M}^c_{O_n}$,
and we let
\[
R(d\eta)=\sum2^{-n}P_{\mu_n}(X_{O_{n}}\in d\eta).
\]
Note that
(\ref{eqDynkin}) implies
\[
P_{\mu}\bar{H}_{x}^{s}(X_{O_{n}})=
\bar{H}_{x}^{s}(\mu)
\]
for $R\times P_{x,S}$-a.e. $(\mu,s)$. By
Fubini's theorem we deduce that there exists a $P_{x,S}$-null set
$\mathcal{N}$ s.t.
%
%
\begin{equation}\label{eqFubini}
P_{\mu}\bar{H}_{x}^{s}(X_{O_{n}})=
\bar{H}_{x}^{s}(\mu) \qquad\mbox{$\forall n$, for $R$-a.e. $\mu
\in\mathcal{M}_{D}^c$, $\forall s\in\mathcal{N}^{c}$.}
\end{equation}
For $s\in\mathcal{N}^{c}$ and $\mu\in\mathcal{M}^c_D$, we choose $O_n$
containing the support of $\mu$ and define
\[
H^{s}_{x}(\mu)=P_{\mu}\bar{H}^{s}_{x}(X_{O_n}).
\]
$H^s_x(\mu)>0$ since this is true of $\bar{H}^s_x$, but we cannot rule
out $H^s_x(\mu)=\infty$.
We set
$H^{s}_{x}(\mu)$ to some arbitrary positive constant for $s\in
\mathcal
{N}$. The definition of $H_{x}^{s}(\mu)$ is independent of the choice of
$O_n$ since, if $O_k\supset O_n$,
then
\begin{eqnarray*}
P_{\mu}\bar{H}_{x}^{s}(X_{O_{k}}) &=&
P_{\mu}P_{X_{O_{n}}}\bar{H}_{x}^{s}(X_{O_k})
\\
&=&P_{\mu}\bar{H}_{x}^{s}(X_{O_{n}}).
\end{eqnarray*}
The first equality is due to the Markov property. The second equality
is due to (\ref{eqFubini}) and the fact that $P_{\mu}(X_{O_{n}}\in
(\cdot))$ is absolutely continuous with respect to $R$, by Theorem
\ref
{theoabscon}.

Clearly, $H^{s}_{x}(\mu)$ is measurable and by (\ref{eqDynkin}), is a
version of
$\frac{dP_{\mu,S}}{dP_{x,S}}(s)$ for each $\mu\in\mathcal{M}_{D}^c$.
To show that $H^{s}_{x}(\mu)$ satisfies property (c), we need to show that
$H^{s}_{x}(\mu)$ is extended $X$-harmonic for each $s\in\mathcal
{N}^{c}$. Let
$\mu$ and $O$ be s.t. $\mu\in\mathcal{M}_{O}^c$ and pick $O_n$ s.t. $O$
is compactly contained in $O_n$.
Then, by definition,
\[
H^{s}_{x}(\mu)=P_{\mu}\bigl(
\bar{H}_{x}^{s}(X_{O_n})\bigr)
\]
and
\[
P_{\mu}H^{s}_{x}(X_{O})=P_{\mu}P_{X_{O}}
\bigl(\bar{H}^{s}_{x}(X_{O_n})\bigr).
\]
By the Markov property these two are equal. Now the family $\{
H_{x}^{s},s\in S\}$ satisfies properties (a), (b) and (c).\vadjust{\goodbreak}

Let us\vspace*{1pt} define $P^{s}_{\mu}$ as in (\ref{eqhtrans}). It remains to
prove that $P^{s}_{\mu}$ is the desired transition kernel. Let
$D'\Subset D$ and
$Y\in\mathcal{F}_{\subset D'}$. Then
\begin{eqnarray*}
P_{\mu}\bigl(f(S)P_{\mu}^S(Y)\bigr) &=&\int
f(s)P_\mu^s(Y) P_{\mu,S}(ds)
\\[-2pt]
&=&\int f(s)\frac{1}{H_x^s(\mu)}\int Y(\omega)H_x^s
\bigl(X_{D'}(\omega)\bigr)P_\mu(d\omega)P_{\mu,S}(ds)
\\[-2pt]
&=&\int Y(\omega)\int\frac{1}{H_x^s(\mu)}f(s)H_x^s
\bigl(X_{D'}(\omega)\bigr)P_{\mu,S}(ds)P_\mu(d
\omega)
\\[-2pt]
&=&\int Y(\omega)\int f(s)H_x^s\bigl(X_{D'}(
\omega)\bigr)P_{x,S}(ds)P_\mu(d\omega)
\\[-2pt]
&=&\int Y(\omega)\int f(s)P_{X_{D'}(\omega),S}(ds)P_\mu(d\omega)
\\[-2pt]
&=&P_\mu\bigl(YP_{X_{D'}}\bigl(f(S)\bigr)\bigr)
\\[-2pt]
&=&P_\mu\bigl(f(S)Y\bigr).
\end{eqnarray*}
Here we are using the definition of $P_{\mu,S}$, the definition of
$P_\mu^s$, Fubini's theorem, the definition of $H_x^s$, the definition
of $P_{\nu,S}$, and the Markov property of $X$.

Uniqueness follows by a similar argument. Suppose $\{H_x^s\}_{s\in
\Sigma}$ and $\{\tilde{H}_x^{s}\}_{s\in\Sigma}$ be any two families
with the properties (a), (b), (c) and (d). Then
\[
H_x^s(\mu)=\tilde H_x^s(\mu)
\qquad\mbox{for $P_{x,S}$-a.e. $s$, $\forall\mu$}
\]
because of property (b).
With $R(d\mu)$ as before, there is therefore a $P_{x,S}$-null set
$\mathcal{N}$ such that
\[
H_x^s(\mu)=\tilde H_x^s(\mu)
\qquad\mbox{for $R$-a.e. $\mu$ and for $s\notin\mathcal{N}$.}
\]
Let $\mu\in\mathcal{M}^c_D$ and $s\notin\mathcal{N}$. Choose $O_n$ such
that $\mu\in\mathcal{M}^c_{O_n}$. Then by absolute continuity and the
property (c),
\[
H_x^s(\mu)=P_\mu H_x^s(X_{O_n})=P_\mu
\tilde H_x^s(X_{O_n})=\tilde
H_x^s(\mu).
\]
\upqed
\end{pf}

In the remainder of Section~\ref{secboundstat}, we consider the three
special cases described above. Our goal is to obtain relatively
explicit formulas for $H_{x}^{s}$ in each case.\vspace*{-2pt}

\subsection{\texorpdfstring{Conditioning on a Poisson random measure with characteristic measure $\beta X_D$}
{Conditioning on a Poisson random measure with characteristic measure beta X D}}
\label{secpoisson}
Let $N$ be a Poisson random
variable with mean $\langle X_D,\beta\rangle$. Let $Z=\{
Z_1,Z_2,\ldots
\}$ be an i.i.d. sequence of random variables from $X_{D}/\langle
X_D,1\rangle$. Let
\[
Y_{\beta}=\sum_{i=1}^{N}
\delta_{Z_i}.
\]
Note that conditioned on $X_D$, $Y_{\beta}$ is a Poisson random measure
with characteristic measure $\beta X_{D}$, and that the construction makes
sense even if $X_D=0$, because then both $N$ and $Y$ equal 0.\vadjust{\goodbreak}

Taking $S=Y_\beta$, Theorem~\ref{theocondprob} gives an extended
$X$-harmonic function (which we denote $H^{\beta,\nu}_x$ to make
explicit the dependence on $\beta$) for conditioning on $Y_\beta=\nu$.
Here $\nu$ is an atomic measure. We let $P^{\beta,\nu}_\mu$ denote the
law of the corresponding conditioned process. In principle this is only
uniquely defined for a.e. $\nu$, but we will find an explicit form that
is valid more generally.

It will be convenient to also define variants of these objects. For any
positive integer $k$, take $S_k=(Z_1,\ldots,Z_k)$ if $N=k$, and
$S_k=\Delta\notin\partial D$ otherwise. Set
$X_{D}^{k}(dz_1,\ldots,dz_k)$ to be the product measure\vspace*{1pt}
$X_{D}(dz_1)\times\cdots\times X_{D}(dz_k)$, and let
$P_{\mu,S_k}^{\beta}$ be the distribution of $S_k$ with respect to
$P_{\mu}$. In other words, for $k\geq1$ and $f\dvtx(\partial D)^k\cup\{
\Delta\}\rightarrow\mathbb{R}$ such that $f(\Delta)=0$, we have
\begin{eqnarray*}
P_{\mu,S_k}^{\beta}(f) &=& P_\mu\bigl(f(S_k)
\bigr)
\\
&=& \frac{1}{k!}P_{\mu} \biggl(\biggl(\int_{(\partial
D)^{k}}
\frac{f}{\langle X_D,1\rangle^k}\,dX_{D}^{k}\biggr)\langle X_D,
\beta\rangle^ke^{-\langle X_{D},\beta\rangle}1_{\{X_D\neq0\}} \biggr)
\\
&=& \frac{\beta^k}{k!}P_{\mu} \biggl(\biggl(\int_{(\partial
D)^{k}}f\,dX_{D}^{k}
\biggr)e^{-\langle X_{D},\beta\rangle} \biggr).
\end{eqnarray*}

So for any $\beta>0$, positive integer $k$, and any $k$-tuple
$z=\{z_i\}$ of elements of~$\partial D$, Theorem~\ref{theocondprob}
gives us a family of extended $X$-harmonic functions
$\{H^{\beta,k,z}\dvtx z=(z_1,\ldots,z_k)\in(\partial D)^{k}\}$ such that
%
%
\begin{equation}\label{kfoldRNdensity}
H^{\beta,k,z}_{x}(\mu)=\frac{dP_{\mu,S_k}^{\beta
}}{dP_{x,S_k}^{\beta
}}(z_1,z_2,\ldots,z_k).
\end{equation}

Let $l_\beta\doteq4V_D(\beta)$. We have $P_{\mu,Y_\beta}\{0\}
=P_\mu
(e^{-\langle X_D,\beta\rangle})
=e^{-\langle\mu,l_\beta\rangle}$ and $P_{x,Y_\beta}\{0\}
=e^{-l_\beta
(x)}$, and therefore we find $H^{\beta,\nu}$ for $\nu=0$ simply by
the ratio
%
%
\begin{equation}\label{eqHat0}
H^{\beta,0}_x(\mu)=\frac{e^{-\langle\mu,l_\beta\rangle
}}{e^{-l_\beta
(x)}}
\end{equation}
by Theorem~\ref{theocondprob}.

Let $l\ge0$ be a bounded Borel function on $D$. For $x\in D$, we let
$m_x^l(dz)=\Pi_x^l(\xi_{\tau_D}\in dz)$ denote harmonic measure on
$\partial D$ for the\vspace*{1pt} operator $L^l$. Then $m_x^l$ and $m_y^l$ are
mutually absolutely continuous, for $x,y\in D$. [This is a well-known
fact; however, for the curious reader, here is a quick argument for why
it is true. Let $D'$ be a smooth domain, compactly contained in $D$,
and $x,y\in D'$. Let $m_{x,D'}^{l}$ be the harmonic measure on
$\partial
{D'}$. If $A$ is a Borel subset of $\partial D$ and $m^l_x(A)=0$,
because of the strong Markov property and the fact that $m_{x,D'}^{l}$
is equivalent to the surface measure $\gamma_{D'}$ on $\partial D'$, we
have that $m_{z}^{l}(A)=0$ for $\gamma_{D'}$ almost all $z$. This
implies $m_{y}^{l}(A)=0$, again due to the strong Markov property and
the fact that $m_{y,D'}^{l}\sim\gamma_{D'}$.]

Let
\[
k^l_x(y,z)=\frac{dm_y^l}{dm_x^l}(z)
\]
denote the density. If $D$ were sufficiently regular, this would be a
version of the Martin kernel for the operator $L^l$, but we make no
such regularity assumptions at this point. We take $k^l$ to be a
jointly measurable version of this density that is harmonic in $y$, for
each $z\in\partial D$. One can construct $k^l$ in a similar way as in
Theorem~\ref{theocondprob}. That is, we start with a family $\{\tilde
{k}((\cdot),z),z\in\partial D\}$ such that $\tilde{k}$ is measurable as a
function of $(y,z)$, and for fixed $y$, $\tilde{k}(y,(\cdot))$ is a version
of $\frac{dm_y^l}{dm_x^l}(z)$. The existence of such a family follows
from Theorem~A.1 of \citet{Dyn06} and the absolute continuity of the
family $\{m_y^l,y\in D\}$ with respect to $m_x^l$. Then we take a
sequence $D_n\Subset D$ exhausting $D$, and let $k^{l}(y,z)=\Pi
_y^l(\tilde{k}(\xi_{\tau_{D_n}},z))$ for $y\in D_n$. Then we prove that
$k^{l}((\cdot),z)$ is well defined, and harmonic for all $z$ except on an
$m_x^l$-null set $\mathcal{N}$, on which we set $k^l$ to be an
arbitrary constant. We omit the details as the arguments are very
similar to those in the proof of Theorem~\ref{theocondprob}.

In the case $l=0$ we write $m_x(dz)=m_x^0(dz)$ and $k_x(y,z)=k_x^0(y,z)$.

The particular case of interest is $l=l_\beta= 4V_D(\beta)$. Suppose
that $k\ge1$ and that $z_1,\ldots,z_k\in\partial D$. For $C\subset
K=\{
1,\ldots,k\}$, recursively define
\[
\rho_{C}^\beta= \cases{ k_{x}^{l_\beta}(
\cdot,z_i),&\quad for $C=\{i\}$,
\vspace*{2pt}\cr
\displaystyle \frac{1}{2}\sum
_{A \subset C, \varnothing\neq A\neq C} G_{D}^{l_\beta
}\bigl(4
\rho_{A}^\beta\rho_{C\setminus A}^\beta\bigr), &\quad
for $|C|>1$.}
\]
Finally, set
\[
\rho_{\mu}^{\beta,k}(z_1,\ldots,z_k)=e^{-\langle\mu,l_\beta
\rangle}
\sum\bigl\langle\mu,\rho^{\beta}_{C_1}\bigr\rangle
\cdots\bigl\langle\mu,\rho^{\beta
}_{C_r}\bigr\rangle,
\]
where the sum ranges over all partitions $\{C_1,\ldots,C_r\}$ of $K$.

In the following theorem, we use the convention that
$H^{\beta,0,z}=H^{\beta,0}$.
%
%
\begin{theorem}
\label{thmPoissonrepn}
Let $D\Subset E$, and $\beta\ge0$ and $x\in D$. Then:

\begin{longlist}[(a)]
\item[(a)]$H_{x}^{\beta,\nu}=H_{x}^{\beta,k,z}$ for $P_{x, Y_{\beta
}}$-almost all finite atomic measures $\nu$, where $k$ and $z$ are
such that
%
%
\begin{equation}\label{eqequiv}
\nu(dx)=\sum_{1}^{k}
\delta_{z_{i}}(dx).
\end{equation}

\item[(b)] For $(m_{x}^{l_\beta})^k$-a.e. $(z_1,\ldots,z_k)$, for all $\mu
\in
\mathcal{M}_D^c$,
%
%
\begin{equation}\label{eqPoisson}
H^{\beta,k,z}_x(\mu)
=\frac{\rho_{\mu}^{\beta,k}(z_1,\ldots,z_k)}{\rho_{x}^{\beta,k}
(z_1,\ldots,z_k)}.
\end{equation}

\item[(c)] If $D$ is smooth, then in fact
$\rho^{\beta,k}_\mu(z_1,\ldots,z_k)<\infty$ for all
$\mu\in\mathcal{M}_D^c$ whenever $z_1,\ldots,z_k$ are distinct.
\end{longlist}
\end{theorem}
\begin{pf}
(a) $P^{\beta}_{\mu,S_k}(f)$ remains unchanged if we permute the
arguments of $f$. Thus we can choose the densities $H^{\beta,k,z}_x(\mu
)$ to be both $X$-harmonic and invariant under permutations of the
$z_i$. A simple way to confirm this is to replace an $X$-harmonic
choice of $H^{\beta,k,z}_x(\mu)$ by $\frac{1}{k!}\sum_\sigma
H^{\beta,k,\sigma(z)}_x(\mu)$, where the sum is over permutations
$\sigma$. The latter is still $X$-harmonic, and a version of the
density $dP^{\beta}_{\mu,S_k}/dP^{\beta}_{x,S_k}$, but is also clearly
invariant under permutations.

For a finite atomic measure $\nu$, all of whose atoms have mass 1, find
$k$ and $z_1,\ldots,z_k$ such that (\ref{eqequiv}) holds. Then define
\[
\tilde{H}^{\beta,\nu}_{x}(\mu):=H_{x}^{\beta,k,z}(
\mu).
\]
Note that $\tilde{H}^{\beta,\nu}_x(\mu)$ is well defined, since
$H_{x}^{\beta,k,z}$ depends
only on $z^{k}:=(z_1,\ldots,\break z_k)$ and is invariant under permuting
$z^{k}$. [Note, if two sequences $z$ and $\tilde{z}$ satisfy
(\ref{eqequiv}), then $z^{k}$ and $\tilde{z}^{k}$ must be
permutations of each other.]

Let $f$ be a function defined on the space of finite atomic measures.
If $\nu=\sum_{i=1}^k\delta_{z_i}$, write $f_k(z)$ for $f(\nu)$. To
finish the proof it is enough to
observe
\begin{eqnarray*}
&&
P_{x, Y_{\beta}}\bigl(\tilde{H}^{\beta,(\cdot)}_{x}(\mu)f(
\cdot)\bigr)
\\
&&\qquad= P_{x}\bigl(\tilde{H}^{\beta,Y_\beta}_{x}(
\mu)f(Y_\beta)\bigr)
\\
&&\qquad=P_{x,Y_\beta}\{0\}\tilde{H}^{\beta,0}(\mu)f(0)
\\
&&\qquad\quad{}+\sum_{k=1}^{\infty} P_{x}
\biggl(e^{-\langle X_D,\beta\rangle}\frac{\langle X_D,\beta
\rangle^k}{k!}1_{\{X_D\neq0\}}\int
H_{x}^{\beta,k,z}(\mu)f_k(z)\frac
{X^k_D(dz)}{\langle X_D,1\rangle^k}
\biggr)
\\
&&\qquad=P_{\mu,Y_\beta}\{0\}f(0)+\sum_{k=1}^{\infty}
\frac{\beta^k}{k!} P_{x} \biggl(e^{-\langle X_D,\beta\rangle} \int
H_{x}^{\beta,k,z}(\mu)f_k(z)X^k_D(dz)
\biggr)
\\
&&\qquad=P_{\mu,Y_\beta}\{0\}f(0)+\sum_{k=1}^{\infty}
P^{\beta}_{x,S_k}\bigl(H_{x}^{\beta,k,(\cdot)}(
\mu)f_k(\cdot)\bigr)
\\
&&\qquad=P_{\mu,Y_\beta}\{0\}f(0)+\sum_{k=1}^{\infty}
P^{\beta}_{\mu,S_k}(f_k)
\\
&&\qquad=P_{\mu,Y_\beta}\{0\}f(0)+\sum_{k=1}^{\infty}
P_{\mu} \biggl(e^{-\langle X_D,\beta\rangle}\frac{\langle X_D,\beta
\rangle^k}{k!}1_{\{X_D\neq0\}} \int
f_k(z)\frac{X^k_D(dz)}{\langle X_D,1\rangle^k} \biggr)
\\
&&\qquad= P_{\mu}\bigl(f(Y_\beta)\bigr)
\\
&&\qquad= P_{\mu, Y_{\beta}}(f).
\end{eqnarray*}

(b) Define $\tilde H^{\beta,k,z}_x(\mu)$ to be the right-hand side of
(\ref{eqPoisson}). Following an argument of \citet{Dyn04a},
Chapter 5,\vspace*{-1pt}
one can show that $\rho^{\beta,k}_{\mu}$ is the density of
$P_{\mu,S_k}^{\beta}$ with respect to $(m_{x}^{l_\beta})^k$. The
argument uses the moment formulas (\ref{eqpalm1}), (\ref{eqrecur1}),
(\ref{eqrecur2}) and then pulls $k$ factors of harmonic measure out of
the resulting expressions, leaving the densities $k_{x}^{l_\beta}$
behind. It follows that $\tilde H^{\beta,k,z}_x(\mu)$ is a version of
the Radon--Nikodym derivative in (\ref{kfoldRNdensity}). The finiteness
condition for $\rho^{\beta,k}_\mu$ follows immediately.

Furthermore, 
by Theorem 3.1 of \citet{SV}, $\tilde H^{\beta,k,z}_x$ is $X$-harmonic;
see remark (iv) below. Thus\vspace*{1pt} $\tilde H^{\beta,k,z}_x= H^{\beta,k,z}_x$
for $(m_{x}^{l_\beta})^k$-a.e. $z$, which is the sense up to which
$H^{\beta,k,z}_x$ is well defined.

(c) The argument for (c) is a straightforward
modification of the estimates used in Theorem 5.3 of \citet{SV}.
\end{pf}

\begin{Remarks*}
(i) The conclusion is that we have obtained an explicit formula for
$H^{\beta,\nu}_x(\mu)$. The abstract definition of this $X$-harmonic
function was valid only up to an unspecified null set of $\nu$'s,
whereas the canonical expression we have obtained is well defined as
long as $\nu$ is a finite atomic measure, all of whose atoms have mass
1 (assuming that $D$ is smooth).

\mbox{}\hphantom{i}(ii) The arguments of this section would work equally well
for conditioning on the value of a Poisson random measure with
characteristic measure $\beta(x)X_D(dx)$, where $\beta(x)$ is now a
bounded measurable function on $\partial D$.

(iii) If $D$ is smooth, then instead of taking $k^l_x(y,z)$
to be the density of $m_y^l$ with respect to $m_x^l$, we could use the
Poisson kernel in its place, and get a similar result. In other words,
we could take the density of $m_y^l$ with respect to the surface
measure $\gamma$ on $\partial D$, rather than the density with respect
to $m_x^l$.

\mbox{}\hspace*{0.6pt}(iv) $\tilde{H}^{\beta,k,z}$ falls in the family of
$X$-harmonic functions
considered in
\citet{SV}. This family of
$X$-harmonic functions are characterized by a function $g$, and
$\mathcal{L}^{4g}$-harmonic
functions $v_1,\ldots,v_k$. In our example the function
$g$ is $u_\beta=V_D\beta$, and the harmonic
functions $v_i$ are the functions $k^{l_\beta}_x(\cdot,z_i)$. In
\citet{SV} it is shown that for $D$ Lipschitz of dimension $d\ge4$,
$g=0$ and $v_i=k_x(\cdot,z_i)$ where $z_1,\ldots,z_k$ are distinct points
chosen on the boundary, the resulting $X$-harmonic function
corresponds to conditioning SBM to hit the points $z_i$. The same
argument would work in dimension $d=3$, at least when $D$ is smooth.
\end{Remarks*}

\subsection{\texorpdfstring{Conditioning on a r.v. $Z$ sampled from measure $\frac{X_D}{{\langle X_D,1\rangle}}$}
{Conditioning on a r.v. $Z$ sampled from measure $(X_D)/(<X_D,1>)$}} \label{secpointsampling}

Recall that the random variable $Z$ is drawn from the probability distribution
$\frac{X_{D}}{{\langle X_D,1\rangle}}$ if $X_D\neq0$, and set equal
to some given
$\Delta
\notin\partial D$ if $X_D=0$. Applying Theorem~\ref{theocondprob}
gives us a family of extended
$X$-harmonic functions
\[
H^{z}_x=\frac{dP_{\mu,Z}}{dP_{x,Z}}(z)
\]
indexed by points $z$ of
$\{\Delta\}\cup\partial D$.
We denote the law of the corresponding conditional process by $P^z_\mu$.

Recall\vspace*{1pt} that $\xi_t$ is a Brownian motion under $\Pi_y$. For $z\in
\partial D$, we let $\Pi_y^z$ be a probability under which $\xi_t$ is a
$k_x(\cdot, z)$-transform of Brownian motion. [Recall $k_x((\cdot
),z)\doteq k^{0}_{x}((\cdot),z)$.] In other words,
\[
\Pi_y^z\bigl(f(\xi_t), t<
\tau_D\bigr)=\frac{1}{k_x(y,z)}\Pi_y\bigl(f(
\xi_t)k_x(\xi_t,z), t<\tau_D
\bigr)
\]
for every bounded measurable $f$.

The following result establishes a concrete formula for $H^z_x$ that is
defined for $m_x$-a.e. $z\in\partial D$ when $D$ is a general domain.
When $D$ is smooth, the same argument as in the previous section gives
a canonical version, defined for all $z\in\partial D$.
%
%
\begin{theorem}
Let $D\Subset E$ and $x\in D$.
Then for
$m_x(dz)$-almost all $z\in\partial D$, $H_x^z(\mu)<\infty$, and
%
%
\begin{equation}
H_{x}^{z}(\mu)=\frac{\int_{0}^{\infty}\langle\mu,k_x(\cdot,z)\Pi_{(\cdot
)}^{z}(e^{-\phi(u_\beta)})\rangle e^{-\langle\mu,u_\beta
\rangle}\,d\beta}{\int_{0}^{\infty}
\Pi_{x}^{z}(e^{-\phi(u_\beta)})e^{-u_\beta(x)}\,d\beta}
\end{equation}
for every $\mu$, where $u_\beta=V_D\beta$ and the random variable
$\phi
(u_\beta)$ is defined as
%
%
\begin{equation}
\phi(u_\beta)=4\int_{0}^{\tau_D}u_{\beta}(
\xi_t)\,dt.
\end{equation}
\end{theorem}
\begin{pf}
We first find the Radon--Nikodym derivative of $P_{\mu,Z}$ w.r.t.
the harmonic measure $m_x(dz)$ on the boundary of $D$. We
observe that
%
%
\begin{eqnarray}\label{eqintder}
P_{\mu,Z}(f)&=& P_{\mu} \biggl(\frac{\langle X_D,f\rangle}{\langle
X_D,1\rangle}1_{\{X_D
\neq
0\}}
\biggr)\nonumber\\[-8pt]\\[-8pt]
&=&-\int_{0}^{\infty}\frac{d}{d\lambda}
P_{\mu} \bigl(e^{-\lambda\langle X_D,f\rangle-\beta
\langle
X_{D},1\rangle
} \bigr)\bigg|_{\lambda=0}
\,d\beta.\nonumber
\end{eqnarray}
Note that the above derivative equals 0 when $X_D=0$.
By the branching property,
%
%
\begin{equation}\label{eqbranch}
P_{\mu}\bigl(e^{-\lambda\langle X_{D},f\rangle-\beta\langle X_D,1\rangle
}\bigr)=e^{-\langle\mu,u_{\lambda
f+\beta}\rangle},
\end{equation}
where
\[
u_{\lambda f+\beta}=V_D(\lambda f+\beta)=\mathbb{N}_{(\cdot
)}
\bigl(1-e^{-\langle X_D, \lambda f+\beta\rangle}\bigr).
\]
Taking the derivative of the right-hand side of (\ref{eqbranch}), and
evaluating at $\lambda=0$ we get
%
%
\begin{equation}\qquad
P_{\mu} \biggl(\frac{\langle X_D,f\rangle}{\langle
X_D,1\rangle}1_{\{X_D \neq0\}} \biggr)=\int_0^\infty\bigl\langle\mu,
\mathbb{N}_{(\cdot)} \bigl(\langle X_{D},f\rangle e^{-\beta\langle
X_{D},1\rangle}\bigr)\bigr\rangle e^{-\langle
\mu,u_\beta\rangle}\,d\beta.
\end{equation}
Differentiation under the integral sign is easily justified. By the
Palm formula,
\begin{eqnarray*}
\mathbb{N}_{y}\bigl(\langle X_{D},f\rangle
e^{-\beta\langle X_{D},1\rangle}\bigr) &=&\Pi_y\bigl(f(\xi_{\tau
_D})e^{-\phi(u_\beta)}
\bigr) \\
&=&
\int_{\partial D}\Pi_{y}^{z}
\bigl(e^{-\phi(u_{\beta})}\bigr)f(z) m_y(dz)
\\
&=&\int_{\partial D}\Pi_{y}^{z}
\bigl(e^{-\phi(u_{\beta})}\bigr)k_x(y,z)f(z) m_x(dz).
\end{eqnarray*}
So,
\[
\bigl\langle\mu, \mathbb{N}_{(\cdot)}\bigl(\langle X_{D},f
\rangle e^{-\beta
\langle X_{D},1\rangle}\bigr)\bigr\rangle=\int_{\partial D}\bigl
\langle\mu,k_x(\cdot,z)\Pi_{(\cdot
)}^{z}
\bigl(e^{-\phi
(u_{\beta})}\bigr)\bigr\rangle f(z) m_x(dz).
\]
Hence
\begin{eqnarray*}
&&
\int_{0}^{\infty}\bigl\langle\mu, \mathbb{N}_{(\cdot)}\bigl(\langle
X_{D},f\rangle e^{-\beta\langle X_{D},1\rangle} \bigr)\bigr\rangle
e^{-\langle\mu,u_{\beta}\rangle}\,d\beta
\\
&&\qquad =\int_{\partial D} f(z) \biggl(\int_{0}^{\infty}
\bigl\langle\mu,k_x(\cdot,z)\Pi_{(\cdot
)}^{z}
\bigl(e^{-\phi(u_\beta)}\bigr)\bigr\rangle e^{-\langle\mu,
u_{\beta}\rangle}\,d\beta\biggr)
m_x(dz).
\end{eqnarray*}

Therefore, both $P_{\mu,Z}$ and $P_{x,Z}$ have densities with
respect to $m_x(dz)$, given by
\[
\int_{0}^{\infty}\bigl\langle\mu,k_x(
\cdot,z)\Pi_{(\cdot
)}^{z}\bigl(e^{-\phi
(u_{\beta})}\bigr)\bigr\rangle
e^{-\langle\mu,u_{\beta}\rangle}\,d\beta
\]
and
\[
\int_{0}^{\infty}\Pi_{x}^{z}
\bigl(e^{-\phi(u_\beta)}\bigr)e^{-u_\beta
(x)}\,d\beta,
\]
respectively. The ratio of these two is a version of the desired
Radon--Nikodym derivative. To show that it equals $H_x^z$ for almost
all $z$, it simply remains to show that it is extended $X$-harmonic.

The denominator is simply a normalizing factor, so consider the
numerator. For $l_\beta=4u_\beta$, it is known [see Theorem 1.1 of
\citet{SV}] that $\mu\mapsto\langle\mu,v\rangle e^{-\langle\mu
,u_{\beta
}\rangle}$ is $X$-harmonic whenever $v$ is $\mathcal{L}^{l_\beta
}$-harmonic. And in our case,
$\langle\mu,k_x(\cdot,z)\Pi_{(\cdot)}^{z}(e^{-\phi(u_{\beta
})})\rangle
e^{-\langle\mu,u_{\beta}\rangle}
=\langle\mu,k_x^{l_\beta}(\cdot,z)\rangle e^{-\langle\mu,u_{\beta
}\rangle}$, as required.

Now we show that $H_{x}^{z}(\mu)<\infty$ for all $\mu$, for
$m_x$-almost all $z$. Let $\mu_0\in\mathcal{M}$ be fixed. Since
$H_{x}^{z}(\mu_0)$ is a density, there exist a $m_x$ null set
$B\subset
\partial D$ s.t. for $z\in B$, $H_{x}^{z}(\mu_0)$ is finite all $z\in
B^c$. Let $\mu\in\mathcal{M}_{D}^c$ and choose $D'\Subset D$ such that
both $\mu$ and $\mu_0$ are compactly supported in $D'$, and assume that
$D'$ is smooth. Let $p_{\mu}$ be the measure defined on $\partial D'$
by $p_{\mu}(f)=P_{\mu}(e^{-\langle X_{D'},u_{\beta}\rangle}\langle
X_{D'}, f\rangle)$. Then our analysis in Section~\ref{secpoisson}\vadjust{\goodbreak}
gives us that $p_{\mu}$ and $p_{\mu_0}$ are equivalent, and the
Radon--Nikodym density is
\[
\frac{dp_{\mu}}{dp_{\mu_0}}(y)=\frac{\langle\mu, \tilde
{k}^{l_{\beta
}}((\cdot),y)\rangle}{\langle\mu_0, \tilde{k}^{l_{\beta}}((\cdot
),y)\rangle},
\]
where $\tilde{k}^{l_{\beta}}(u,y)$ is the Poisson kernel of $D'$ for
the operator $\mathcal{L}^{l_{\beta}}$. Note that this density is
bounded by a constant $C(\mu,\mu_0)$ since $\tilde{k}^{l_{\beta}}(u,y)$
is harmonic in $u$ on the support of $\mu$ and $\mu_0$. Hence
\begin{eqnarray*}
\bigl\langle\mu,k_x^{l_\beta}(\cdot,z)\bigr\rangle
e^{-\langle\mu,u_{\beta
}\rangle}&=&P_{\mu}\bigl\langle X_{D'},k_x^{l_\beta}(
\cdot,z)\bigr\rangle e^{-\langle X_{D'},u_{\beta}\rangle}
\\
&\leq& C(\mu,\mu_0) P_{\mu_0}\bigl\langle
X_{D'},k_x^{l_\beta}(\cdot,z)\bigr\rangle
e^{-\langle X_{D'},u_{\beta}\rangle}
\\
&=& C(\mu,\mu_0)\bigl\langle\mu_0,k_x^{l_\beta}(
\cdot,z)\bigr\rangle e^{-\langle
\mu,u_{\beta}\rangle}.
\end{eqnarray*}
It follows that $H_{x}^{z}(\mu)\leq C(\mu,\mu_0)H_{x}^{z}(\mu_0)<\infty
$ for all $z\in B$, and hence the proof is complete.
\end{pf}

\subsection{Conditioning on a linear function of $X_{D}$}\label{seclinear}

Let $L$ be a linear and measurable map from the linear cone of positive
finite measures
$\mathcal{M}_{\partial D}$ on $\partial D$ to a Luzin measurable space
$(V,\mathcal{V})$ where $V$ is a vector space, and $\mathcal{V}$ is
countably generated.
Let $V_+$ be the image of $\mathcal{M}_{\partial D}$, and write
$V^*=V_+\setminus\{0\}$. Assume that $L\mu=0$ implies $\mu=0$.

Let $T_{n}$ be the map $V^{n}_+\to V^{n}_+$ defined by
\[
T_{n}(v_1,\ldots,v_n)\mapsto
(v_1+v_2+\cdots+v_n,v_1,\ldots,v_{n-1}),
\]
and for
$A\in\mathcal{V}$ and $A\subset V^*$, let
\[
N_{x,L(X_{D})}(A)=\mathbb{N}_{x}\bigl(L(X_{D})\in A,
X_D\neq0\bigr).
\]
We fix $x\in D$ as before and define a reference measure $R_x$ on $V^*$ by
\[
R_x(A)=P_{x,L(X_{D})}(A, X_D\neq0).
\]
Its total mass is $r_{x,0}=P_{x,L(X_D)}(V^*)=1-e^{-u(x)}$, where
$u(x)=\break-\log P_x(X_D=0)$. Note that $u=\lim_{\beta\rightarrow\infty}
V_D(\beta)$.
Throughout this section we will assume that
%
%
\begin{equation}\label{eqregdom}\quad
\mbox{$D$ is a bounded domain, all of whose boundary points are regular.}
\end{equation}
This holds, for example, if the boundary of $D$ is smooth. Under
assumption~(\ref{eqregdom}), $V_D(\beta)$ is the unique solution of
\[
\tfrac12\Delta u=2u^2
\]
on $D$ with $u=\beta$ on $\partial{D}$ [Proposition 8.2.1.B of
\citet{Dyn02}]. Because of this and the comparison principle
[Proposition 8.2.1.H of \citet{Dyn02}], $V_D(\beta)$ is also the
maximal solution of $\frac12\Delta u=2u^2$ on $D$ bounded by $\beta$.
We will need this in the proof of Lemma~\ref{theoref}.

By Theorem~\ref{theocondprob} we know the existence of a family of
extended $X$-harmonic functions $\{H^v_x,v\in V\}$ such that $H^v_x(\mu
)=dP_{\mu,L(X_D)}/dP_{x,L(X_D)}(v)$. In this section we are going to
find a more explicit formula for this family.
%
%
\begin{lemma}\label{theoref} Assume (\ref{eqregdom}).
There exists a family of functions $\{\gamma_{x,v}\dvtx D\mapsto
(0,\infty),
v\in V^*\}$ such that the mapping $(v,y)\mapsto\gamma_{x,v}(y)$ is
measurable and for all $y\in D$
%
%
\begin{equation}\label{eqgamma}
N_{y,L(X_{D})} (dv)=\gamma_{x,v}(y)R_x(dv).
\end{equation}
In addition, there exists
a measurable kernel $K_{x,n}(v;dv_1,dv_2,\ldots,dv_{n-1})$
from $V^*$ to $(V^*)^{n-1}$, such that
%
%
\begin{eqnarray}\label{eqkernel}
&&
R_{x}^{n}\circ T_{n}^{-1}(dv,
dv_1,dv_2,\ldots,dv_{n-1})\nonumber\\[-8pt]\\[-8pt]
&&\qquad=K_{x,n}(v;dv_1,dv_2,
\ldots,dv_{n-1})R(dv).\nonumber
\end{eqnarray}
Moreover $K_{x,n}(v,\cdot)$ is a strictly positive measure, for
$R$-a.e. $v$.
\end{lemma}
\begin{pf} Recall that
\[
P_{\mu,L(X_{D})}(A)=P_{\mu}\bigl(L(X_{D})\in A\bigr).
\]
Because all $P_{\mu}(X_D\in\cdot)$ are equivalent, so are the
$P_{\mu,L(X_{D})}$, as are their restriction to $V^*$. Thus all
$P_{\mu,L(X_{D})}$ (when restricted to $V^*$) are equivalent to $R_x$.
Moreover, since
\[
N_{y,L(X_D)}(A)=\mathbb{N}_{y}\bigl(L(X_{D})\in
A,X_D\neq0\bigr)=\mathbb{N}_{y}\bigl(P_{X_{D'}}
\bigl(L(X_{D})\in A,X_D\neq0\bigr)\bigr)
\]
for all $y\in D$ and $D'\Subset D$ such that $y\in D'$,
$N_{y,L(X_{D})}$ is also equivalent to
$R_x$ for all $y\in D$.
By Theorem A.1 of \citet{Dyn06} we get a family of functions $\{
\gamma
_{x,v}\dvtx D\mapsto[0,\infty), v\in V^*\}$ such that the mapping
$(v,y)\mapsto\gamma_{x,v}(y)$ is measurable, and (\ref{eqgamma})
holds. Clearly such a $\gamma_{x,v}(y)$ can be chosen strictly positive
since $N_{y,L(X_{D})}$ and $R_x$ are equivalent.

It will be convenient for the proof to write $R_x(dv)=r_{x,0}\tilde
R_x(dv)$, where $\tilde R_x$ is a probability measure. In other words,
$\tilde R_x(dv)=P_{x}(L(X_D)\in dv\mid X_D\neq0)$. Note that with this
choice of $\tilde R_x$, $\tilde R_{x}^{n}\circ
T_{n}^{-1}(dv, dv_1,dv_2,\ldots,dv_{n-1})$
is the joint distribution of
$(V_1+\cdots+V_{n}, V_1,\ldots,V_{n-1})$ where $V_{i}$ are
independent random variables with distribution $\tilde R_x$. Let
$R_{x,n}$ be
the marginal distribution of $V_1+\cdots+V_{n}$, where the $V_i$ are as
above. The following
decomposition is then immediate:
\[
\tilde R_{x}^{n}\circ T_{n}^{-1}(dv,
dv_1,dv_2,\ldots,dv_{n-1})=\tilde
{K}_{x,n}(v;dv_1,dv_2,\ldots,dv_{n-1})R_{x,n}(dv),
\]
where $\tilde{K}_{x,n}(v;dv_1,dv_2,\ldots,dv_{n-1})$ is the
conditional probability kernel for $(V_1,\ldots,V_{n-1})$ given
$V_1+\cdots+V_{n}$.

We now show that $R_{x,n}$ is absolutely continuous with respect to $R_x$.
Let $X^{1}_{D},\ldots,X^{n}_{D}$ be $n$
independent realizations of
the exit measure under the law $P_{x}$. Then the distribution of
$X^{1}_{D}+\cdots+{X}^{n}_{D}$ is given by the $P_{n\delta_x}$
distribution of
$X_{D}$.

Let $F$ be s.t. $R_{x}(F)=0$, that is,
%
%
\begin{equation}\label{eqR1}
P_{x}\bigl(F\bigl(L(X_{D})\bigr)1_{\{X_{D}\neq0\}}
\bigr)=0.
\end{equation}
Because $P_x$ and $P_{n\delta_x}$ are absolutely continuous, (\ref{eqR1})
implies
%
%
\begin{equation}\label{eqR2}
P_{n\delta_x}\bigl(F\bigl(L(X_{D})\bigr)1_{\{X_{D}\neq0\}}
\bigr)=0.
\end{equation}
Since
\begin{eqnarray*}
&&
P_{n\delta_x}\bigl(F\bigl(L(X_{D})\bigr)1_{\{X_{D}\neq
0\}}
\bigr)
\\
&&\qquad=(P_{x})^n\bigl(F\bigl(L\bigl(X^{1}_{D}+
\cdots+X^{n}_{D}\bigr)\bigr)1_{\{X^{1}_{D}+\cdots+X^{n}_{D}\neq
0\}}\bigr)
\\
&&\qquad\geq(P_x)^{n}\bigl(F\bigl(L\bigl(X^{1}_{D}
\bigr)+\cdots+L\bigl(X^{n}_{D}\bigr)\bigr)1_{\{X^{1}_{D}\neq0\}}
\cdots1_{\{X^{n}_{D}\neq0\}}\bigr)
\\
&&\qquad=r_{x,0}^{n}R_{x,n}(F),
\end{eqnarray*}
this implies that $R_{x,n}(F)=0$, so indeed, $R_{x,n}$ is absolutely
continuous with respect to $R$.

If $h_{x}^{n}(v)$
is the Radon--Nikodym derivative of $R_{x,n}$ with
respect to $R_x$, we get (\ref{eqkernel}) with
\[
K_{x,n}(v;dv_1,dv_2,\ldots,dv_{n-1})=
\tilde{K}_{x,n}(v;dv_1,dv_2,\ldots,dv_{n-1})h_{x}^{n}(v).
\]

It remains only to show that $K_{x,n}(v,\cdot)$ is strictly positive.
Because\break $\tilde K_{x,n}(v,\cdot)$ is, this amounts to showing the
converse to the absolute continuity result above, namely that $R_x$ is
absolutely continuous with respect to~$R_{x,n}$.

Our approach is to use the Poisson representation, as in the absolute
continuity argument in \citet{Dyn04a}. Suppose $R_{x,n}(F)=0$ and
$0\le
F\le1$. The Poisson representation gives that
\begin{eqnarray*}
&&
P_\mu\bigl(F\bigl(L(X_D)\bigr), X_D
\neq0\bigr)
\\
&&\qquad=\sum_{k=1}^\infty\frac{e^{-\langle\mu,u\rangle}}{k!}\iint F\bigl(L(\nu_1)+\cdots+L(\nu_k)\bigr)
\N_{x_1}(X_D\in d\nu_1, X_D\neq0)\cdots
\\
&&\hspace*{67pt}\qquad\quad{} \times\N_{x_n}(X_D\in d\nu_k,
X_D\neq0)\mu(dx_1)\cdots\mu(dx_k)
\\
&&\qquad=\sum_{k=1}^\infty\frac{e^{-\langle\mu,u\rangle}}{k!}\int
f_k(x_1,\ldots,x_k)\mu(dx_1)
\cdots\mu(dx_k),
\end{eqnarray*}
where
\begin{eqnarray*}
f_k(x_1,\ldots,x_k)&=&\int F\bigl(L(
\nu_1)+\cdots+L(\nu_k)\bigr)
\N_{x_1}(X_D\in d\nu_1, X_D
\neq0) \cdots\\
&&\hspace*{10pt}{}\times \N_{x_k}(X_D\in d\nu_k,
X_D\neq0).
\end{eqnarray*}
Let $D_m\Subset D$ such that $x\in D_m$ and $D_m\uparrow D$. Then
\begin{eqnarray*}
&&
P_x\bigl(F\bigl(L(X_D)\bigr), X_D
\neq0\bigr)
\\
&&\qquad=P_x\bigl(P_{X_{D_m}}\bigl(F\bigl(L(X_D)
\bigr), X_D\neq0\bigr)\bigr)
\\
&&\qquad=P_x\Biggl(\sum_{k=1}^\infty
\frac{e^{-\langle X_{D_m},u\rangle
}}{k!}\int f_k(x_1,\ldots,x_k)X_{D_m}(dx_1)
\cdots X_{D_m}(dx_k)\Biggr).
\end{eqnarray*}
There is a similar Poisson representation for $R_{x,n}(F)$, involving a
sum of integrals of the $f_k$ for $k\ge n$. Since $R_{x,n}(F)=0$, we
conclude that for each $k\ge n$ there are $x_1,\ldots,x_k\in D$ such
that $f_k(x_1,\ldots,x_k)=0$. By absolute continuity, we conclude that
$f_k(x_1,\ldots,x_k)=0$ for every $x_1,\ldots, x_k$.

Since $F\le1$ we obtain the bound
\[
f(x_1,\ldots,x_k)\le\prod_{j=1}^k
\N_{x_j}(X_D\neq0)=u(x_1)\cdots
u(x_k).
\]
Therefore
\[
P_x\bigl(F\bigl(L(X_D)\bigr), X_D\neq0
\bigr) \le\sum_{k=1}^{n-1}P_x
\biggl(\frac{e^{-\langle X_{D_m},u\rangle
}}{k!}\langle X_{D_m},u\rangle^k\biggr).
\]
The result will follow once we argue that all terms $e^{-\langle
X_{D_m},u\rangle}\langle X_{D_m},u\rangle^k$ converge to 0 $P_x$-a.s.
as we let $m\rightarrow\infty$, by the dominated convergence theorem
since these terms are bounded. To show that $e^{-\langle
X_{D_m},u\rangle}\langle X_{D_m},u\rangle^k\rightarrow0$, it is enough
to show that the stochastic boundary value of $u$ (i.e., $\lim
_{m\rightarrow\infty} \langle X_{D_m},\break u\rangle$) is $0$ or $\infty$,
$P_x$-a.s. Let $Z_{\beta}$ be the stochastic boundary value of the
constant function $\beta$. The sequence $\langle X_{D_m},\beta\rangle$
is a uniformly integrable martingale with respect to $P_y$ for all
$y\in D$. That it is a martingale follows because the constant function
$\beta$ is harmonic. The uniform integrability follows because this
martingale is square integrable. Indeed, by the moment formula
(\ref{eqrecur2}), $P_y( \langle X_{D_m},\beta\rangle)^{2}=G_{D_m}(4\beta
^{2})(y)+ \beta\leq G_{D}(4\beta^{2})(y)+ \beta<\infty$, since $D$ is
bounded. It follows that the log-potential of $Z$ is the maximal
solution of $\frac12\Delta u=2u^2$ on $D$ bounded by $\beta$, which is
$V_{D}{\beta}$, as we argued at the beginning of this section; see
Sections 9.2.1 and 9.2.2 of \citet{Dyn02}. So $V_{D}\beta(x)=-\log P_x
e^{-Z_{\beta}}$. Also note that $Z_\beta=\beta Z_1$, and $Z_{\beta
}\uparrow Z$ where
%
%
\begin{equation}
Z=\cases{0, &\quad if $Z_1=0$,
\cr
\infty, &\quad if $Z_1>0$.}
\end{equation}

By the dominated convergence theorem
\[
P_x e^{-Z}=\lim_{\beta
\rightarrow
\infty}P_x e^{-Z_{\beta}}=\lim_{\beta\rightarrow\infty}
e^{-V_D\beta(x)}=e^{-u(x)},
\]
so $u$ is the log potential of $Z$, and
therefore $Z$ is the stochastic boundary value of~$u$. Since $Z$ is
$0$ or $\infty$ $P_x$-a.s., the proof is complete.
\end{pf}
%
%
\begin{lemma}\label{theokerdecompose} Assume (\ref{eqregdom}). Let
$(K_{x,n})_{n\geq2}$ be a sequence of transition kernels satisfying
(\ref{eqkernel}) for $n\geq2$. Define
%
%
\begin{equation}
\bar{K}_{x,n}(v;dv_1,\ldots,dv_n):=K_{n}(v;dv_1,
\ldots,dv_{n-1})\times\delta_{v-(v_1+\cdots+v_{n-1})}(dv_n).\hspace*{-25pt}
\end{equation}

Then for $R_x$-almost all $v \in V^*$, the following holds for all
$n\geq2$, $1\leq r\leq n$, and any partition $C_1,\ldots,C_r$ of
$\{1,\ldots,n\}$, where $n_i=|C_i|$:
\[
\bar{K}_{x,n}(v;dv_1,\ldots,dv_n)=\int
\bar{K}_{x,r}(v;dv_1,\ldots,dv_r)\prod
_{i=1}^{r}\bar{K}_{x,n_i}(\tilde{v}_i;dv_{C_i}).\vspace*{-2pt}
\]
\end{lemma}
\begin{pf} Since $\mathcal{V}$ is countably generated, so is $\mathcal
{V}^n$ (the product $\sigma$-field), and therefore for each $n\geq1$,
there exists a sequence of nonnegative Borel measurable functions $\{
f^{n}_{i}\}_{i=1}^{\infty}$ generating $\mathcal{V}^n$.

It will suffice to show for any $k,j\geq1$, $n\geq2$, $1\leq r\leq
n$, and any given partition $C_1,\ldots,C_r$ of $\{1,\ldots,n\}$ that
\begin{eqnarray*}
&&
\int f^{1}_{j}(v) \biggl[\int
f^n_k(v_1,\ldots,v_n)\bar
{K}_{x,n}(v,dv_1,\ldots,dv_n)
\biggr]R_x(dv)
\\[-2pt]
&&\qquad=\int f^1_j(v) \Biggl[\int\bar{K}_{x,r}(v,d
\tilde{v}_1,\ldots,d\tilde{v}_r)
\\[-2pt]
&&\hspace*{52.5pt}\qquad\quad{}\times\int\prod_{i=1}^{r}
\bar{K}_{x,n_i}(\tilde{v}_i,dv_{C_i})
f_k^n(v_1,\ldots,v_n)
\Biggr]R_x(dv).
\end{eqnarray*}
Let
\[
\tilde{f}^n_k(\tilde{v}_1,\ldots,
\tilde{v}_r)=\int\prod_{i=1}^{r} \bar
{K}_{x,n_i}(\tilde{v}_i,dv_{C_i})
f_k^n(v_1,\ldots,v_n).
\]
Then
\begin{eqnarray*}
&&
\int f^1_j(v)\int\bar{K}_{x,r}(v,d
\tilde{v}_1,\ldots,d\tilde{v}_r)\int\prod
_{i=1}^{r} \bar{K}_{x,n_i}(
\tilde{v}_i,dv_{C_i}) f_k^n(v_1,\ldots,v_n)R_x(dv)
\\[-2pt]
&&\qquad=\int f^{1}_{j}(\tilde{v}_1+\cdots+
\tilde{v}_{r})\tilde{f}_{k}^{n}(
\tilde{v}_1,\ldots,\tilde{v}_r)R_x^r(d
\tilde{v}_1,\ldots,d\tilde{v}_r)
\\[-2pt]
&&\qquad=\int f^1_j(\tilde{v}_1+\cdots+
\tilde{v}_{r})f_k^n(v_1,\ldots,v_n) \Biggl[\prod_{i=1}^{r}
\bar{K}_{x,n_i}(\tilde{v}_i,dv_{C_i})R_x(d
\tilde{v}_i) \Biggr]
\\[-2pt]
&&\qquad=\int f^1_j(\Sigma_{C_1}v_m+
\cdots+ \Sigma_{C_r}v_m)f^n_k
(v_1,\ldots,v_n) \prod_{i=1}^{r}
\prod_{C_i}R_x(dv_m)
\\[-2pt]
&&\qquad=\int f^{1}_{j}(v_1+\cdots+v_n)f^n_k(v_1,\ldots
,v_n)R_x^n(dv_1,\ldots,dv_n)
\\[-2pt]
&&\qquad=\int f^{1}_{j}(v) \biggl[\int f^n_k(v_1,\ldots,v_n)\bar
{K}_{x,n}(v,dv_1,\ldots,dv_n) \biggr]R_x(dv).\hspace*{50pt}\qed
\end{eqnarray*}
\noqed
\end{pf}
%
%
\begin{lemma} \label{theonullset1} Assume (\ref{eqregdom}). For any
$B_0\subset V^*$ with $R_x(B_0^c)=0$, there exists $B\subset V^*$ with
$R_x(B^c)=0$, such that $B\subset B_0$ and for all $n\geq2$ and $v\in
B$, we will have $(v_1,\ldots,v_n) \in B^{n}$, $\bar
{K}_{x,n}(v,dv_1,\ldots,dv_n)$ a.s.
\end{lemma}
\begin{pf}
Define recursively $B_{m}, m\geq1$ as
\[
B_{m}=\Biggl\{v\in B_{m-1}\dvtx\sum
_{n=2}^{\infty}\int\sum_{i=1}^{n}1_{B_{m-1}^c}(v_i)
\bar{K}_{x,n} (v,dv_1,\ldots,dv_n)= 0\Biggr\}.
\]
Then $R_x(B_{m}^c)=0$. Because $R_{x}(B_0^c)=0$, and assuming
$R_x(B_{m-1}^c)=0$, we have that
\begin{eqnarray*}
&&
\int\sum_{n=2}^{\infty}\int\sum
_{i=1}^{n}1_{B_{m-1}^c}(v_i)
\bar{K}_{x,n} (v,dv_1,\ldots, dv_n)R_x(dv)
\\[-2pt]
&&\qquad=\sum_{n=2}^{\infty}\int\sum
_{i=1}^{n}1_{B_{m-1}^c}(v_i)
R_x^n(dv_1,\ldots,dv_n)
\\[-2pt]
&&\qquad=\sum_{n=2}^{\infty}nR_x
\bigl(B_{m-1}^c\bigr)r_{x,0}^{n-1}=0,
\end{eqnarray*}
which implies $R_x\{v\dvtx\sum_{n=2}^{\infty}\int\sum
_{i=1}^{n}1_{B_{m-1}^c}(v_i)\bar{K}_{x,n} (v,dv_1,\ldots, dv_n)\neq
0\}
=0$. This implies that $R_x(B_{m}^c)=0$ since
\[
B_{m}^{c}\subset\Biggl\{v\dvtx\int\sum
_{i=1}^{n}1_{B_{m-1}^c}(v_i)
\bar{K}_{x,n} (v,dv_1,\ldots, dv_n)\neq0\Biggr
\}.
\]
Let $B=\bigcap_{m=0}^{\infty} B_{m}$. Clearly $R_x(B^c)=0$ and
$B\subset B_0$.

Now, if $v\in B$, then for any $m$, $\sum_{n=2}^{\infty}\int
\sum_{i=1}^{n}1_{B_{m}^c}(v_i)\bar{K}_{x,n} (v,dv_1,\ldots, dv_n) =0$.
By the monotone convergence theorem, this implies that
\[
\sum_{n=2}^{\infty}\int\sum
_{i=1}^{n} 1_{B^c}(v_i)
\bar{K}_{x,n} (v,dv_1,\ldots,dv_n)=0.
\]
Hence for all $n$,
$\bar{K}_{x,n}(v,dv_1,\ldots,dv_n)$ almost all $(v_1,\ldots,v_n)$ is
in $B$.
\end{pf}

So far $\gamma_{x,v}(y)$ is any jointly measurable version of the
Radon--Nikodym derivative of $N_{y,X_D}$ with respect to $R_x$. In the
following theorem, we refine this choice and find a formula for the
family of extended $X$ harmonic functions $\{H_{x}^{v},v\in V^*\}$
corresponding to $dP_{\mu,L(X_D)}/dP_{x,L(X_D)}(v)$.
%
%
\begin{theorem} \label{theolinear} Assume (\ref{eqregdom}). There
exists a version of $\{\gamma_{x,v}, v\in V^*\}$ of Lemma~\ref
{theoref} such that for $R_x$ almost every $v\in V^*$
%
%
\begin{equation}\label{eqversion}
\gamma_{x,v}(y)=\N_{y}\bigl(H_{x}^{v}(X_{D'})
\bigr)
\end{equation}
for every $y$, and $D'\Subset D$ such that $y\in D'$, where
%
%
\begin{equation}
\label{eqlinear} H_{x}^{v}(\mu)=\cases{
\displaystyle e^{-\langle\mu,u\rangle+u(x)},
\hspace*{116pt} \qquad\mbox{if $v=0$},
\vspace*{2pt}\cr
\displaystyle e^{-\langle\mu,u\rangle}\langle\mu,\gamma_{x,v}\rangle\cr
\displaystyle \qquad{}+ \sum
_{n=2}^{\infty}\int\frac{1}{n!}e^{-\langle\mu,u\rangle
}K_{x,n}(v;dv_1,
\ldots,dv_{n-1})
\vspace*{2pt}\cr
\qquad\hspace*{39pt}\displaystyle {}\times
\langle\mu,\gamma_{v_1}\rangle\cdots\langle{\mu,\gamma_{x,v_{n-1}}}%
\rangle\vspace*{2pt}\cr
\qquad\hspace*{39pt}\displaystyle {}\times\langle{\mu,
\gamma_{x,v-(v_1+\cdots+v_{n-1})}}\rangle, \qquad \mbox{if
$v\neq0$}.}
\end{equation}
Moreover, with such a version of $\gamma_{x,v}$, $\{H_{x}^{v},v\in
V^*\}
$ defined by (\ref{eqlinear}) is extended $X$-harmonic for
$R_x$-almost all $v$. For fixed $\mu$, $H_{x}^{v}(\mu)$ is a version of
$dP_{\mu,L(X_D)}/dP_{x,L(X_D)}(v)$.
\end{theorem}
\begin{pf} Let $\{\bar{\gamma}_{x,v},v\in V^*\} $ be a family of
functions satisfying the properties of Lemma~\ref{theoref}. Define
$\bar{H}_{x}^{v}$ by formula (\ref{eqlinear}) with this $\bar
{\gamma
}_{x,v}$ in place of $\gamma_{x,v}$.
Since $\bar{\gamma}_{x,v}>0$, it follows that $\bar{H}_x^v>0$. We
will
first show that for each $\mu$, $\bar{H}_{x}^{v}(\mu)$ is a version of
$dP_{\mu,L(X_D)}/dP_{x,L(X_D)}(v)$.
Let $F$ be a nonnegative Borel function on~$V_+$.
\begin{eqnarray*}
&&
\int P_{x,L(X_D)}(dv)F(v) \bar{H}_{x}^{v}(
\mu)
\\[-2pt]
&&\qquad=P_{x,L(X_D)}\bigl(\{0\}\bigr)F(0)\bar{H}_{x}^{0}(
\mu)+\int F(v) \bar{H}_{x}^{v}(\mu)R_x(dv)
\\[-2pt]
&&\qquad= e^{-{u(x)}}F(0)e^{-\langle\mu,u\rangle
+u(x)}+e^{-\langle\mu,u\rangle}\int F(v)\langle\mu,
\bar{\gamma}_{x,v}\rangle R_x(dv)
\\[-2pt]
&&\qquad\quad{} +\sum_{n=2}^{\infty}\frac{e^{-\langle\mu,u\rangle
}}{n!}\iint\langle\mu,\bar{\gamma}_{x,v_1}\rangle\cdots\langle{\mu,\bar
{\gamma
}_{x,v_{n-1}}} \rangle\langle{\mu,\bar{\gamma}_{x,v-(v_1+\cdots
+v_{n-1})}}\rangle
\\[-2pt]
&&\hspace*{81.2pt}\qquad\quad{} \times F(v)K_{x,n}(v;dv_1,\ldots,dv_{n-1})R_{x}(dv)
\\[-2pt]
&&\qquad=e^{-\langle\mu,u\rangle}F(0)+\sum_{n=1}^{\infty}
\frac
{e^{-\langle\mu,u\rangle}}{n!}\int F(v_1+\cdots+v_n)\langle\mu,\bar{
\gamma}_{x,v_1}\rangle\cdots\langle\mu,\bar{\gamma}_{x,v_n}
\rangle
\\[-2pt]
&&\hspace*{128pt}\qquad\quad{} \times R_x(dv_1)\cdots R_x(dv_n)
\\[-2pt]
&&\qquad=e^{-\langle\mu,u\rangle}F(0)\\[-2pt]
&&\qquad\quad{}+\sum_{n=1}^{\infty}
\frac
{e^{-\langle\mu,u\rangle}}{n!}\int F(v_1+\cdots+v_n)
\\[-2pt]
&&\hspace*{75pt}\qquad\quad{} \times
\bigl\langle\mu,N_{(\cdot),L(X_{D})}(dv_1)\bigr
\rangle\cdots
\bigl\langle\mu,N_{(\cdot),L(X_{D})} (dv_n)\bigr\rangle
\\[-2pt]
&&\qquad=\sum_{n=0}^{\infty} \frac{e^{-\langle\mu,u\rangle}}{n!}\int
F\bigl(L(v_1)+\cdots+L(\nu_n)\bigr)\mathcal{R}_\mu(d
\nu_1)\cdots\mathcal{R}_\mu(d\nu_n)
\\[-2pt]
&&\qquad=P_\mu\bigl(F\bigl(L(X_D)\bigr)\bigr) =P_{\mu,L(X_{D})}(F),
\end{eqnarray*}
where we are using the Poisson representation.
We choose countable base $O_n$ (w.l.o.g. closed under finite unions)
and consider the measure $R$ on $\mathcal{M}_{D}^{c}$ defined as in the
proof of Theorem~\ref{theocondprob}.
The argument in that proof tells us that
there exists a $P_{x,L(X_D)}$-null set
$\mathcal{N}$ s.t.
%
%
\begin{equation}\label{eqFubini2}
P_{\mu}\bar{H}_{x}^{v}(X_{O_{n}})=
\bar{H}_{x}^{v}(\mu) \qquad\mbox{$\forall n$, for $R$-a.e. $\mu
\in\mathcal{M}_{D}^c$, $\forall v\in\mathcal{N}^{c}$.}
\end{equation}
Without loss of generality we may assume that for $v\in\mathcal
{N}^{c}$, the statement of Lemma~\ref{theokerdecompose} holds, and
moreover by Lemma~\ref{theonullset1}, $\bar{K}_n(v,dv_1,\ldots, dv_n)$
almost surely $(v_1,\ldots,v_n)\in(\mathcal{N}^{c})^n$ for all $n$. For
$v\in\mathcal{N}^{c}$ and $y \in D$, we choose $O_n$ containing the
support of $y$ and define
\[
\gamma_{x,v}(y)=\N_y\bigl(\bar{H}^{v}_{x}(X_{O_n})
\bigr).
\]
Since $\bar{H}_x^v>0$ it follows that $\gamma_{x,v}(y)\in(0,\infty]$
for every $y$.
We set
$\gamma_{x,v}$ to some arbitrary constant $>0$ for $v\in\mathcal{N}$.
The definition of $\gamma_{x,v}$ is independent of the choice of
$O_n$ since, if $O_k\supset O_n$,
then
\begin{eqnarray*}
\N_y\bar{H}_{x}^{v}(X_{O_{k}}) &=&
\N_y P_{X_{O_{n}}}\bar{H}_{x}^{s}(X_{O_k})
\\[-2pt]
&=&\N_y \bar{H}_{x}^{s}(X_{O_{n}}).
\end{eqnarray*}
The first equality is due to Markov property. The second equality
is due to (\ref{eqFubini2}) and the fact that $\N_y(X_{O_{n}}\in
(\cdot))$ is absolutely continuous with respect to $R$.

We now show that for fixed $y$, this gives a version of $\frac
{dN_{y,L(X_D)}}{dR_x}$.
If $y\in O_n$, then
\begin{eqnarray*}
\int_{\{v\neq0\}} F(v) N_{y,L(X_D)}(dv) &=&\N_y
\bigl(F\bigl(L(X_D)\bigr)1_{\{X_D\neq0\}} \bigr)
\\[-2pt]
&=&\N_y \bigl(P_{X_{O_n}}\bigl(F\bigl(L(X_D)
\bigr)1_{\{X_D\neq0\}}\bigr) \bigr)
\\[-2pt]
&=&\N_y\biggl(\int_{\{v\neq0\}}\bar{H}_{x}^{v}(X_{O_n})F(v)R_x(dv)
\biggr)
\\[-2pt]
&=&\int_{\{v\neq0\}}F(v)R_x(dv)\N_{y}
\bigl(H_{x}^{v}(X_{O_n})\bigr)
\\[-2pt]
&=&\int_{\{v\neq0\}}F(v)\gamma_{x,v}(y)R_x(dv).
\end{eqnarray*}
Let $H^{v}_{x}$ be defined by formula (\ref{eqlinear}). Then we know
that for fixed $\mu$, $H^{v}_{x}(\mu)$ is a version of $dP_{\mu
,L(X_D)}/dP_{x,L(X_D)}(v)$. Now we\vspace*{1pt} are going to show that for $v\in
\mathcal{N}^{c}$, $H^{v}_{x}(\mu)=P_{\mu}(\bar{H}_{x}^v(X_{O_k}))$. To
simplify the notation, we drop the basepoint $x$ from $\bar{K}_{x,n}$,
$\bar{\gamma}_{x,v}$ and $\gamma_{x,v}$ for the remainder of the
proof,
\begin{eqnarray*}
\hspace*{-4pt}&&
P_\mu\bar{H}_{x}^{v}(X_{O_k})
\\[-2pt]
\hspace*{-4pt}&&\quad=P_\mu\sum_{n=1}^{\infty}\int
\frac{e^{-\langle
X_{O_k},u\rangle}}{n!} \prod_{i=1}^{n}\langle X_{O_k},
\bar{\gamma}_{v_i}\rangle\bar{K}_{n}(v;dv_1,
\ldots,dv_n)
\\[-2pt]
\hspace*{-4pt}&&\quad=\sum_{n=1}^{\infty}
\frac{1}{n!}\int P_\mu\Biggl(e^{-\langle
X_{O_k},u\rangle}
\prod_{i=1}^{n}\langle X_{O_k},\bar{\gamma
}_{v_i}\rangle\Biggr) \bar{K}_{n}(v;dv_1,
\ldots,dv_n)
\\[-2pt]
\hspace*{-4pt}&&\quad=\sum_{n=1}^{\infty}\frac{1}{n!}\int
e^{-\langle\mu,V_{O_k}(u)\rangle} \Biggl[\sum_{\pi(n)} \prod
_{i=1}^{r}\biggl\langle\mu,\mathbb{N}_{(\cdot)}
\biggl(e^{-\langle
X_{O_k},u\rangle}\prod_{j\in C_i}\langle X_{O_k},\bar{
\gamma}_{v_{j}}\rangle\biggr)\biggr\rangle\Biggr]
\\[-2pt]
\hspace*{-4pt}&&\hspace*{38pt}\qquad{}\times\bar{K}_{n}(v;dv_1,
\ldots,dv_n)
\\[-2pt]
\hspace*{-4pt}&&\quad=e^{-\langle\mu,V_{O_k}(u)\rangle}\\[-2pt]
\hspace*{-4pt}&&\qquad{}\times\sum_{n=1}^{\infty}
\frac{1}{n!}\sum_{\pi(n)} \int\bar
{K}_{r}(v;d\tilde{v}_1,\ldots,\tilde{v}_r)
\\[-2pt]
\hspace*{-4pt}&&\hspace*{57pt}\qquad{}\times\prod_{i=1}^{r}
\biggl\langle\mu,\mathbb{N}_{(\cdot)} \biggl(e^{-\langle
X_{O_n},u\rangle} \biggl[\int
\prod_{j\in C_i}\langle X_{O_n},\bar{\gamma}_{v_{j}}\rangle
\bar{K}_{n_i}(\tilde{v}_i;dv_{C_i}) \biggr]
\biggr)\biggr\rangle
\\[-2pt]
\hspace*{-4pt}&&\quad=e^{-\langle\mu,V_{O_k}(u\rangle)
}\sum_{n=1}^{\infty}
\frac{1}{n!}\sum_{r=1}^{n}\int
\bar{K}_{r}(v;\tilde{v}_1,\ldots,\tilde{v}_r)
\sum_{n_1,\ldots,n_r}\frac
{n!}{n_1!\cdots
n_r!r!}
\\[-2pt]
\hspace*{-4pt}&&\qquad{}\times\prod_{i=1}^{r}\biggl
\langle\mu,\mathbb{N}_{(\cdot)} \biggl(e^{-\langle
X_{O_k},u\rangle} \biggl[\int\prod
_{j\in C_i}\langle X_{O_k},\bar{
\gamma}_{v_{j}}\rangle\bar{K}_{n_i}(\tilde{v}_i;dv_{C_i})
\biggr] \biggr)\biggr\rangle
\\[-2pt]
\hspace*{-4pt}&&\quad=e^{-\langle\mu,V_{O_k}(u\rangle)}\sum_{r=1}^{\infty}
\int\frac{1}{r!}\bar{K}_{r}(v;d\tilde{v}_1,\ldots,d
\tilde{v}_r)
\\[-2pt]
\hspace*{-4pt}&&\qquad{}\times\prod_{i=1}^{r}\sum
_{n_i=1}^{\infty}\frac
{1}{n_i!}\biggl\langle
\mu,\mathbb{N}_{(\cdot)} \biggl(e^{-\langle
X_{O_k},u\rangle} \biggl[\int\prod
_{j\in C_i}\langle X_{O_k},\bar{\gamma}_{v_{j}}
\rangle\bar{K}_{n_i}(\tilde{v}_i;dv_{C_i})
\biggr] \biggr)\biggr\rangle
\\[-2pt]
\hspace*{-4pt}&&\quad=e^{-\langle\mu,V_{O_k}(u\rangle)}\sum_{r=1}^{\infty}
\int\frac{1}{r!}\bar{K}_{r}(v;d\tilde{v}_1,\ldots,d
\tilde{v}_r)
\\[-2pt]
\hspace*{-4pt}&&\qquad{}\times\prod_{i=1}^{r}\Biggl
\langle\mu,\mathbb{N}_{(\cdot)} \Biggl(\sum_{n_i=1}^{\infty}
\frac{1}{n_i!}e^{-\langle
X_{O_k},u\rangle} \biggl[\int\prod_{j\in C_i}
\langle X_{O_k},\bar{\gamma}_{v_{j}}\rangle\bar{K}_{n_i}(
\tilde{v}_i;dv_{C_i}) \biggr] \Biggr)\Biggr\rangle
\\[-2pt]
\hspace*{-4pt}&&\quad=e^{-\langle\mu,u\rangle}\sum_{r=1}^{\infty}
\int\frac{1}{r!}\bar{K}_{r}(v;d\tilde{v}_1,\ldots,d
\tilde{v}_r)\prod_{i=1}^{r}
\langle\mu,\gamma_{x,\tilde{v}_i}\rangle
\\[-2pt]
\hspace*{-4pt}&&\quad=H_{x}^{v}(\mu),
\end{eqnarray*}
where the last line follows because by our assumption that if $v\in
\mathcal{N}^c$, then $\bar{K}_r(v, d\tilde{v}_1,\ldots,d\tilde{v}_r)$
almost surely $(\tilde{v}_1,\ldots,\tilde{v}_r)\in(\mathcal{N}^c)^n$
and by construction for each $\tilde{v}_i$,
\[
\gamma_{\tilde{v}_i}(y)= \mathbb{N}_{y}\sum
_{n_i=1}^{\infty}\frac{1}{n_i!}e^{-\langle
X_{O_k},u\rangle}\int
\prod_{j\in C_i}\langle X_{O_k},\bar{\gamma}_{v_{j}}\rangle
\bar{K}_{n_i}(\tilde{v}_i;dv_{C_i}).
\]

Now we show that $\gamma_{x,v}$ satisfies equation (\ref{eqversion}).
Let $y\in\tilde{D}\Subset D$. Then for some~$n$, $\tilde{D}\Subset
O_n$. By the definition of $\gamma_{x,v}$ and the Markov property,
\begin{eqnarray*}
\N_y\bigl(H^{v}_{x}(X_{\tilde{D}})
\bigr)&=&\N_y\bigl(P_{X_{\tilde{D}}}\bar{H}_{x,v}(X_{O_n})
\bigr)\\
&=&\N_y \bigl(\bar{H}_{x,v}(X_{O_n})\bigr)=
\gamma_{x,v}(y)
\end{eqnarray*}
as desired.

Replacing $\mathbb{N}_y$ with $P_{\mu}$ in the above argument we have
that $H_{x}^{v}$ is extended $X$-harmonic for all $v\in\mathcal{N}^c$.
It is also clear that $H_{x}^{v}$ is $X$-harmonic for $v=0$. Hence
the proof is complete.
\end{pf}

\begin{Remarks*}
(i) Although Theorem~\ref{theolinear} gives a workable form for
$H^v_x(\mu)$, it does not remove all ambiguity in the
choice of $H^v_x$, since $\gamma_{x,v}$ and
$K_x(v;\cdot)$ are only well defined for a.e. $v$. Ideally we would
like to prove continuity properties of these objects in $v$ as well,
that would then specify them uniquely. But we have not succeeded in
doing that. In subsequent sections we will, however, be able to clarify
the structure of these objects, and show how they determine the
behavior of the $H^v_x$-transformed super-Brownian motion. We do have
regularity in $y$. In particular, we will soon see that the version of
$\gamma_{x,v}(y)$ given by Theorem~\ref{theolinear} is already lower
semi continuous in $y$.

\mbox{}\hphantom{i}(ii) Note that in Sections~\ref{secpoisson} and~\ref{secpointsampling}
without much difficulty we were able to show that the corresponding
extended $X$-harmonic functions are finite, therefore $X$-harmonic.
Finiteness is harder to prove for the extended $X$-harmonic functions
of this section because it requires analytical bounds on the
Radon--Nikodym densities of $n$-order moment measures of SBM for all
$n\geq1$. We hope to pursue such bounds in a subsequent paper.

(iii) An important case is when
$L(X_{D})=X_{D}$. The corresponding extended $X$-harmonic function
$H_{x}^{\nu}$ can be thought as the analogue of the Martin kernel.

\mbox{}\hspace*{0.5pt}(iv) A more tractable application should be the
case $L(X_D)=\langle X_D,1\rangle$, where we condition on the total
mass. In that case, $\gamma_{x,v}(y)$ is a function of only
finite-dimensional variables $v\in[0,\infty)$ and $y\in D$. We hope to
explore this example further in a subsequent paper.

\mbox{}\hspace*{0.5pt}\hphantom{i}(v) An interesting direction is to
explore the relationship between $H_{x}^{\beta,\nu}$of Section
\ref{secpoisson} and $H_{x}^{\nu}$. In particular, with $\beta=n$, what
happens to $H^{n,\nu_{n}}$ as $n\rightarrow\infty$, if
$\{\nu_n\}_{n\geq1}$ is a sequence of finite atomic measures such that
$n^{-1}\nu_{n}$ converges to a finite measure $\nu$ on the boundary?
One can show that
\[
P_{\mu}^{n,Y_n}\rightarrow P_\mu^{X_{D}}
\]
weakly almost surely, in a sense, and use this to investigate whether
$H_{x}^{\nu}$ is extreme. We will pursue this direction in a
subsequent paper.
\end{Remarks*}

\section{\texorpdfstring{Fragmentation system description of $P_{\mu}^{\nu}$}
{Fragmentation system description of P mu nu}}\label{setotalmass}

The results\vspace*{1pt} of this section apply in general to the
conditional law $P_{\mu}^{v}$ given\vadjust{\goodbreak} $L(X_{D})=v$. For simplicity,
however, we will carry out the computations for the case $L(\nu)=\nu$.
So in this section $V^*$ is the set of positive finite measures on
$\partial{D}$.

As in Section~\ref{seclinear}, we assume (\ref{eqregdom}), that is,
that $D$ is regular. From Section~\ref{seclinear}, recall that
$u(x)=-\log P_x(X_D=0)=V_D(\infty)$, and $u$ is a solution of the
boundary value problem
\[
\tfrac12\Delta u=2u^2
\]
on $D$ and $u=\infty$ on $\partial{D}$. We also have
%
%
\begin{equation}
u(x)=\N_x(X_{D}>0)=\lim_{\beta\to\infty}V_D
\beta(x).
\end{equation}

Let $\gamma_{x,\nu}$ and $H_{x}^{\nu}$ be as constructed in Section
\ref{seclinear}. Let $\mathcal{N}\subset V^*$ be the $R_x$-null set
such that for $\nu\in\mathcal{N}^{c}$, $\gamma_{x,\nu}$ and $H_{x}^\nu$
satisfy the system of equations in (\ref{eqversion}) and~(\ref
{eqlinear}), $H_{x}^{\nu}$ is extended $X$-harmonic and the kernels
$\bar{K}_{x,n}$ are strictly positive and satisfy the decomposition
property of Lemma~\ref{theokerdecompose}. By Lemma~\ref{theonullset1}
we may assume that if $\nu\in\mathcal{N}^{c}$, then $\bar
{K}_{x,n}(\nu,d\nu_1,\ldots, d\nu_n)$-almost all $(\nu_1,\ldots,\nu_n)$
are in $(\mathcal{N}^{c})^n$. We fix $\nu\in\mathcal{N}^c$. Let $y$ be
a point in $D$ such that $\gamma_{x,\nu}(y)<\infty$. Recall that
$\gamma_{x,\nu}(y)>0$ by construction. If $y\in D'\Subset D$, we may
then define a change of measure by
\[
\N^\nu_y(Z)=\frac{1}{\gamma_{x,\nu}(y)}\N_y
\bigl(ZH_x^\nu(X_{D'})\bigr)
\]
for positive $\mathcal{F}_{\subset D'}^{y}$-measurable $Z$. ($\mathcal
{F}_{\subset D'}^{y}$ is defined as the $\sigma$-algebra generated by
$\{X_{O}, y\in O \subset D'\}$.)
Since $H_{x}^{\nu}$ is extended $X$-harmonic, $\N_{y}^{\nu}$ is
defined consistently on $\mathcal{F}^{y}_{\subset D-}$, and we have that
$\N_{y}^{\nu}$ is a probability law because of equation (\ref{eqversion}).

In the remainder of this section we turn to
the problem of giving an explicit probabilistic construction of
$\N_y^\nu$ in terms of a backbone along which unconditioned mass
is created. We do this in two steps. Let $H_{x,1}^\nu(\mu
)=e^{-\langle
\mu,u\rangle}\langle\mu,\gamma_{x,\nu}\rangle$, and for $n\ge2$ let
%
%
\begin{equation}
H^{\nu}_{x,n}(\mu)=\int\frac{e^{-\langle\mu,u\rangle
}}{n!}
\bar{K}_{x,n}(\nu;d\nu_1,\ldots,d\nu_{n})
\langle\mu,\gamma_{x,\nu_1}\rangle\cdots\langle{\mu,\gamma_{x,\nu
_{n}}}
\rangle.
\end{equation}
Then $H^{\nu}_{x}=\sum_{n\ge1}H^{\nu}_{x,n}$.
We will first use the recursive moment formula to establish an
inductive relationship for $H^{\nu}_{x,n}$, and then
compare this to the inductive relationship coming from the first
branch of the backbone.

We will make use of a stochastic process $\xi_t$ under various
measures $\Pi_y$ or $\Pi_y^{4u}$. In either case we use the shorthand
\[
\mathcal{N}_{t}^{D'}(\phi)=e^{-\int_0^t 4\N_{\xi_s}(1-e^{-\langle
X_{D'},\phi\rangle}) \,ds},
\]
where $D'\Subset D$ and let $\tau_{D'}$ be the exit
time of $\xi$ from $D'$. We similarly let
\[
\mathcal{N}_{t}^{I,D_I}(\phi_I)=e^{-\int_0^{t\wedge\tau_{D_I}} 4\N
_{\xi
_s}(1-e^{-\langle
X_{I},\phi_I\rangle}) \,ds},
\]
where $I=\{D_1,\ldots,D_k\}$ is such that each $D_j\Subset D$, $\phi
_I=(\phi_1,\ldots,\phi_k)$, $X_I=(X_{D_1},\ldots,X_{D_k})$, so
$\langle
X_I,\phi_I\rangle=\sum_j\langle X_{D_j},\phi_j\rangle$. Also write
$D_I=D_1\cap\cdots\cap D_k$.

In the rest of the paper we will denote $K_{x,2}(\nu,d\nu')$ by
$K_{x}(\nu,d\nu')$. Let $I=\{D_1,\ldots,D_k\}$ where $D_i\subset
D_k=D'\Subset D$. Define a family of operators $N_{y}^{I,\nu,n}$ as
follows:

For $\operatorname{card}(I)=1$,
%
%
\begin{equation}
\label{eqNcardI11} N_y^{I,\nu,n}(\phi)=\cases{
\N_y \bigl(e^{-\langle X_{D'},\phi\rangle} H^{\nu}_{x,1}(X_{D'})
\bigr), &\quad$y\in D'$,
\cr
\gamma_{x,\nu}(y), &\quad$y\notin
D', n=1$,
\cr
0, &\quad$y\notin D',n>1$.}
\end{equation}

For $\operatorname{card}(I)>1$,
%
%
\begin{equation}
\label{eqNcardIk} N_y^{I,\nu,n}(\phi_I)=
\cases{\N_y \bigl(e^{-\langle X_{I},\phi_I\rangle} H^{\nu}_{x,n}(X_{D_k})
\bigr), &\quad$y\in D_I$,
\cr
N_y^{I^j,\nu,n}, &\quad$y
\notin D_j,j\neq k$,}
\end{equation}
where $I_j=I-\{D_j\}$. Note that we have
\begin{eqnarray*}
&&
N_y^{I,\nu,n}(\phi_I)\\
&&\qquad=\frac{1}{n!}\int
n_{C_n}^{I}\bigl(\phi_1,\ldots,
\phi_k+u,\gamma^{\nu_1},\ldots,\gamma^{\nu_n}\bigr)
(y)\\
&&\hspace*{21pt}\qquad\quad{}\times\bar{K}_{x,n}(\nu,d\nu_1,\ldots,d\nu_n),
\end{eqnarray*}
where $C_n=\{1,\ldots,n\}$, and the $n_C^I$ are the operators defined
by equations (\ref{defnmoment1}) and~(\ref{defnmoment2}).
%
%
\begin{lemma}\label{recursion} Assume (\ref{eqregdom}). Let $\phi
_I^u=(\phi_1,\ldots,\phi_{k-1},\phi_k+u)$.
For $n=1$ and \mbox{$y\in D_I$},
%
%
\begin{equation}\label{eqN1}
N_y^{I,\nu,1}(\phi_I)=\Pi_y
\bigl( N_{\xi_{\tau_{D_I}}}^{I,\nu,1}(\phi_I)\mathcal{N}_{\tau_{D_I}}^{I,D_I}
\bigl(\phi_I^u\bigr) \bigr).
\end{equation}
For $n\ge2$ and $y\in D_I$,
%
%
\begin{eqnarray}\label{eqNn}
&&
N_y^{I,\nu,n}(\phi_I)
\nonumber
\\
&&\qquad=\sum_{m=1}^{n-1} \int
\Pi_y \biggl(\int_{0}^{\tau_{D_I}}2
N_{\xi_t}^{I,\nu',m}(\phi_I)N_{\xi_t}^{I,\nu-\nu',n-m}(
\phi_I)\mathcal{N}_{t}^{I,D_I}\bigl(
\phi_I^u\bigr)\,dt \biggr)\nonumber\\[-8pt]\\[-8pt]
&&\hspace*{61pt}{}\times K_x\bigl(\nu,d
\nu'\bigr)
\nonumber\\
&&\qquad\quad{}+\Pi_y N_{\xi_{\tau_{D_I}}}^{I,\nu,n}(
\phi_I)\mathcal{N}_{\tau
_{D_I}}^{I,D_I}\bigl(\phi_I^u\bigr).\nonumber
\end{eqnarray}
\end{lemma}
\begin{pf} If $n=1$, then $H^{\nu}_{x,1}(\mu)=e^{-\langle\mu,u\rangle
}\langle
\mu,\gamma_{x,\nu}\rangle$ and the result is an immediate
consequence of
the basic Palm formula (\ref{eqpalm1}) [for $\operatorname{card}(I)=1$] and the
extended Palm formula (\ref{eqexpalm1}) [for $\operatorname{card}(I)>1$].

If $n\ge2$, by the recursive moment formulas (\ref{eqrecur1}),
\begin{eqnarray*}
\N_y\bigl(e^{-\langle X_{I},\phi_I\rangle} H^{\nu}_{x,n}(X_{D_k})
\bigr)&=&\frac{1}{n!}\int\N_y\bigl(e^{-\langle X_{I},\phi_I\rangle
}e^{-\langle
X_{D_k},u\rangle}
\Pi_i\langle X_{D_k},\gamma_{x,\nu_{i}}\rangle\bigr)\\
&&{}\hspace*{24pt}\times
\bar{K}_{x,n}(\nu;d\nu_1,\ldots,d\nu_{n-1})
\\
&=&A+B,%
\end{eqnarray*}
where
\begin{eqnarray*}
A &=&\frac{1}{2\cdot n!}\int\mathop{\sum_{M \subset N}}_{ \varnothing,
N \ne M}
\Pi_y \biggl(\int_0^{\tau_{D_I}}4
\mathcal{N}_{t}^{I,D_I}\bigl(\phi_I^u
\bigr)
\\
&&\hspace*{112.6pt}{}\times\N_{\xi_t}\bigl(e^{-\langle X_I,\phi_I^u\rangle}\Pi_{i\in
M}
\langle X_{D'},\gamma_{x,\nu_i}\rangle\bigr)
\\
&&\hspace*{112.6pt}{}\times\N_{\xi_t}\bigl(e^{-\langle X_I,\phi_I^u\rangle}\Pi
_{i\notin
M}
\langle X_{D'},\gamma_{x,\nu_i}\rangle\bigr) \,dt \biggr) \\
&&\hspace*{67pt}{}\times\bar
{K}_{x,n}(\nu;d\nu_1,\ldots,d\nu_n).
\end{eqnarray*}
There are ${n\choose m}$ possible choices of $M$ in the above
expression, with cardinality $m$.
Therefore, by rearranging the indices, and using Lemma \ref
{theokerdecompose} with $r=2$, we get
\begin{eqnarray*}
A &=&\sum_{m=1}^{n-1}\int
K_x\bigl(\nu;d\nu'\bigr) \\
&&\hspace*{16.3pt}{}\times\Pi_y \Biggl(
\int_0^{\tau
_{D_I}}2\mathcal{N}_{t}^{I,D_I}
\bigl(\phi_I^u\bigr)
\\
&&\hspace*{56pt}\hspace*{16.3pt}{}\times\Biggl[\frac{1}{m!}\int\N_{\xi_t}
\Biggl(e^{-\langle
X_I,\phi
_I^u\rangle}\prod_{i=1}^m \langle
X_{D_k},\gamma_{x,\nu_i} \rangle\Biggr)\\
&&\hspace*{139.5pt}{}\times \bar{K}_{x,m}
\bigl(\nu';d\nu_1,\ldots,d\nu_{m}\bigr)
\Biggr]
\\
&&\hspace*{56pt}\hspace*{16.3pt}{}\times\Biggl[\frac{1}{(n-m)!}\int\N_{\xi_t}
\Biggl(e^{-\langle
X_I,\phi
_I^u\rangle}\prod_{i=1}^{n-m} \langle
X_{D_k}, \gamma_{x,\nu_i}\rangle\Biggr)
\\
&&\hspace*{69.6pt}
\hspace*{56pt}\hspace*{16.3pt}{}\times\bar{K}_{x,n-m}\bigl(\nu-\nu';d
\nu_1,\ldots,d\nu_{n-m}\bigr) \Biggr]\,dt \Biggr).
\end{eqnarray*}
The term $B$ is $0$ for $\operatorname{card}(I)=1$ and for $\operatorname{card}(I)>1$ is
found using the extended moment formula (\ref{eqexrecur1}) as
\begin{eqnarray*}
B&=& \frac{1}{ n!}\int\Pi_y \bigl( n_{C_n}^{I}(
\phi_1,\ldots,\phi_k+u,\gamma_{x,\nu_1},\ldots,
\gamma_{x,\nu_n},\xi_{\tau
_{D_I}})\mathcal{N}_{\tau_{D_I}}^{D_I}
\bigl(\phi_I^u\bigr) \bigr)
\\
&&\hspace*{21.5pt}{}\times\bar{K}_{x,n}(\nu;d\nu_1,\ldots,d
\nu_{n})
\\
&=&\Pi_y \bigl(N_{\xi_{\tau_{D_I}}}^{I,\nu,n}(
\phi_I)\mathcal{N}_{\tau
_{D_I}}^{I,D_I}\bigl(
\phi_I^u\bigr) \bigr).
\end{eqnarray*}
\upqed
\end{pf}

Now set
$\Gamma_{x,\nu}=2\int\gamma_{x,\nu'}\gamma_{x,\nu-\nu'}
K_x(\nu;d\nu')$.
%
%
\begin{theorem} \label{theogamma} Assume (\ref{eqregdom}). The
function $\gamma_{x,\nu}$ is $\mathcal{L}^{4u}$-superharmonic, and
hence lower-semi-continuous in $y$. For $R_x$-almost all $\nu\in
\mathcal{N}^{c}$ it is in fact an $\mathcal{L}^{4u}$-potential, and satisfies
%
%
\begin{equation}\label{potential}
\gamma_{x,\nu}=G_{D}^{4u} \biggl[2\int
\gamma_{x,\nu}\gamma_{x,\nu
-\nu'} K_x\bigl(\nu;d
\nu'\bigr) \biggr].
\end{equation}
\end{theorem}
\begin{pf} Let $D_k$ be a sequence of domains exhausting $D$. Since
$V_{D_k}u=u$,
we have that $\mathcal{N}_{t}^{D_k}(u)=e^{-\int_0^t4u(\xi_s) \,ds}$ for
$t<\tau_{D_k}$.
Thus
\begin{eqnarray*}
\gamma_{x,\nu}(y)&=&\N_y\bigl(H^{\nu}_{x}(X_{D_k})
\bigr)
=\sum_{n=1}^\infty\N_y
\bigl(H^{\nu}_{x,n}(X_{D_k})\bigr)
\\
&=&\N_{y}\bigl(H_{x,1}^{\nu}(X_{D_k})
\bigr)\\
&&{}+\sum_{n=2}^{\infty} \sum
_{m=1}^{n-1}\int K_x\bigl(\nu; d
\nu'\bigr)
\\
&&\hspace*{46pt}\hspace*{0pt}{}\times
\Pi_y \biggl(\int_0^{\tau_{D_k}}2
\mathcal{N}_{t}^{D_k}(u)\N_{\xi_t}
\bigl(H^{\nu'}_{x,m}\bigr) \N_{\xi_t}
\bigl(H^{\nu-\nu'}_{x,n-m}\bigr) \,dt \biggr)
\\
& = & \N_{y}\bigl(H_{x,1}^{\nu}(X_{D_k})
\bigr) +\sum_{m=1}^{\infty} \sum
_{j=1}^{\infty}\int K_x\bigl(\nu; d
\nu'\bigr)
\\
&&\hspace*{116pt}{}\times\Pi_y^{4u} \biggl(\int
_{0}^{\tau_{D_k}}2\N_{\xi
_t}\bigl(H_{x,m}^{\nu'}
\bigr) \N_{\xi_t}\bigl(H_{x,j}^{\nu-\nu'}\bigr) \,dt \biggr)
\\
& = & \N_{y}\bigl(H_{x,1}^{\nu}(X_{D_k})
\bigr) +2\int\Pi_y^{4u}\biggl(\int_0^{\tau_{D_k}}
\gamma_{x,\nu'}(\xi_t)\gamma_{x,\nu-\nu
'}(
\xi_t) \,dt\biggr) K_x\bigl(\nu; d\nu'
\bigr)
\\
&=&\N_{y}\bigl(H_{x,1}^{\nu}(X_{D_k})
\bigr)+\Pi_y^{4u} \biggl(\int_0^{\tau
_{D_k}}
\Gamma_\nu(\xi_t) \,dt \biggr).
\end{eqnarray*}
The first term is $\mathcal{L}^{4u}$-harmonic on $D_k$ by Lemma \ref
{recursion}, and the second term is an $\mathcal{L}^{4u}$-potential, so
$\gamma_{x,\nu}$ is $\mathcal{L}^{4u}$-superharmonic on each $D_k$.
Thus it is so on $D$ as well.

Moreover,
\begin{eqnarray*}
\int\N_{y}\bigl(H_{x,1}^{\nu}(X_{D_k})
\bigr)R_x(d\nu)&=&\N_{y}%
\biggl(e^{-\langle X_{D_{k}}, u\rangle}\int\langle X_{D_k},%
\gamma_{x,\nu}\rangle R_x(d\nu)\biggr)
\\
&=&\N_{y}%
\biggl(e^{-\langle X_{D_{k}}, u\rangle}\iint
X_{D_k}(dw)%
\gamma_{x,\nu}(w) R_x(d
\nu)\biggr)
\\
&=&\N_{y}%
\biggl(e^{-\langle X_{D_{k}}, u\rangle}\int
X_{D_k}(dw)%
\N_w(X_D\neq0)\biggr)
\\
&=& \N_{y}\bigl(e^{-\langle X_{D_{k}}, u\rangle}\langle X_{D_k},u\rangle
\bigr)
\\
&=&e^{u(y)}P_{y}\bigl(e^{-\langle X_{D_{k}}, u\rangle}\langle
X_{D_k},u\rangle\bigr),
\end{eqnarray*}
where the last equality follows from the recursive moment formula (\ref
{eqrecur2}). As we argued in the proof of Lemma~\ref{theoref}, the
stochastic boundary value of $u$ is 0 or $\infty$ $P_{y}$-a.s.;
therefore $P_{y}(e^{-\langle X_{D_k}, u\rangle}\langle
X_{D_k},u\rangle)$ converge to $0$ as $k\rightarrow\infty$ by the
dominated convergence theorem.
Thus
\[
\int\N_{y}\bigl(H_{1}^{\nu}(X_{D_k})
\bigr)R_x(d\nu)\to0
\]
as $k\to\infty$. By Fatou's lemma,
%
%
\begin{equation}
\label{eqfatou} \liminf_{k\to\infty}\N_{y}\bigl(H_{1}^{\nu}(X_{D_k})
\bigr)=0
\end{equation}
for $R_x$-a.e. $\nu$. But the second term in our expression for
$\gamma_{x,\nu}(y)$ is monotone in $k$, and therefore $\lim_{k\to\infty
}\N_{y}(H_{1}^{\nu}(X_{D_k}))$ must exist and thus equal to $0$ for
$R_x$-a.e. $\nu$ by (\ref{eqfatou}). So, we get that $\gamma_{x,\nu
}(y)=\Pi^{4u}_y(\int_0^{\tau_{D_k}}\Gamma_{x,\nu}(\xi_t) \,dt)$ for
$R_x$-a.e. $\nu$, for every~$y$.

Choosing a countable dense set $y_1, y_2, \ldots\,$, we therefore have a
set $\mathcal{N}_0\supset\mathcal{N}$ such that
$R_x(\mathcal{N}_{0}^{c})=0$, and the above equality holds for every
$y_i$ and for every $\nu\in(\mathcal{N}_0)^c$. Since both functions
are lower-semi-continuous, and agree on a countable dense set, it
follows that (\ref{potential}) holds for every $\nu\in
(\mathcal{N}_0)^c$.
\end{pf}

Suppose now that $\nu\in\mathcal{N}_0^c$ where $\mathcal{N}_0$ is the
$R_x$-null set described in Theorem~\ref{theogamma}. Again we may
assume that if $\nu\in\mathcal{N}_0^c$, then $\bar{K}_{x,n}(\nu,d\nu
_1,\ldots, d\nu_n)$-almost all $(\nu_1,\ldots,\nu_n)$ are in
$(\mathcal
{N}^{c})^n$. Let $\hat\N_y$ be the excursion measure of a super-process
whose spatial motion is killed at rate $u$. In other words, for any
$D'\Subset D$, we have
\[
\hat\N_y\bigl(F(X_{D'})\bigr)=\N_{y}
\bigl(e^{-\langle X_{D'},u\rangle
}F(X_{D'}) \bigr).
\]

Let $y\in D$ be such that $\gamma_{x,\nu}(y)<\infty$. We define a
probability $Q_y^\nu$ on an auxiliary probability space $\tilde
{\Omega
}$ where there is a branching diffusion on $D$, and conditional on this
branching diffusion, a Poisson random measure is generated on the
infinite product space $\mathcal{M}^{\mathcal{O_{D-}}}$. We endow
$\mathcal{M}^{\mathcal{O_{D-}}}$ with the $\sigma$-algebra $\tilde
{\mathcal{F}}_{\subset D-}$ generated by the coordinate maps $\tilde
{x}_{D'}$, $D'\Subset D$ [i.e., $\tilde{x}_{D'}(\omega)=\omega_{D'}$
for $\omega\in\mathcal{M}^{\mathcal{O_{D-}}}$]. Our goal is to
construct an $\mathcal{M}^{\mathcal{O_{D-}}}$-valued process $\tilde
{X}=(\tilde{X}_{D'})_{D'\Subset D}$ on $\tilde{\Omega}$ such that the
law of $\tilde{X}$ with respect to $Q_y^\nu$ will be the same as the
$\mathbb{N}_{y}^{\nu}$ law of the exit measures $(X_{D'})_{D'\Subset
D}$ of a SBM. First we describe how the branching diffusion evolves: we
start a $\gamma_{x,\nu}$-transform of a $\mathcal{L}^{4u}$ process off
at $y$. Since
$\gamma_{x,\nu}$ is a potential, this process dies before reaching
$\partial D$. Say it dies at $w$. Then almost surely $\Gamma_{x,\nu
}(w)$ is finite. Because
\[
\Pi_{y}^{u,\gamma_{x,\nu}}(1_{\Gamma_{x,\nu}(\xi_\zeta)=\infty
})=\frac
{1}{\gamma_{x,\nu}(y)}\int
_{0}^{\infty}\Pi_y^u
\bigl(1_{\Gamma_{x,\nu
}(\xi
_t)=\infty}\Gamma_{x,\nu}(\xi_t)1_{\zeta>t}
\bigr)\,dt,
\]
which\vspace*{1pt} must be equal to $0$ since $\Pi_y^u(1_{\Gamma_{x,\nu}(\xi
_t)=\infty}\Gamma_{x,\nu}(\xi_t)1_{\zeta>t})$ is $0$ or $\infty$\break for
any~$t$, and the assumption\vspace*{1pt} $\gamma_{x,\nu}(y)=\Pi_y^u(\Gamma_{x,\nu
}(\xi_t)1_{\zeta>t})<\infty$ implies\break $\Pi_y^u(1_{\Gamma_{x,\nu
}(\xi
_t)=\infty}\*\Gamma_{x,\nu}(\xi_t)1_{\zeta>t})=0$ for Lebesgue-almost
all $t$, making the right-hand side of this equation equal to $0$.
Note also $\Gamma_{x,\nu}>0$ since $\nu\in\mathcal{N}_0^c$, and
therefore $K_{x}(\nu,d\nu')$ is a strictly positive measure. Now we can
choose $\nu'$ at random with distribution
density
%
%
\begin{equation}\label{eqfragdensity}
\frac{2}{\Gamma_{x,\nu}(w)}\gamma_{x,\nu'}(w)\gamma_{x,\nu-\nu
'}(w)K_{x}
\bigl(\nu;d\nu'\bigr).
\end{equation}

Almost surely, both $\nu'$ and $\nu-\nu'$ are in $\mathcal{N}_0^c$, and
$\gamma_{x,\nu}(w)$ and $\gamma_{x,\nu'}(w)$ are finite, so we may
start two new processes at $w$, following $\gamma_{x,\nu}$ and
$\gamma_{x,\nu-\nu'}$ transforms of~$\mathcal{L}^{4u}$. Note that we
may repeat this process infinitely often.
This defines a branching particle system. Let $\Upsilon_t$ denote
the measure-valued process putting a unit point mass at the historical
paths of
each particle alive at time $t$. (A historical path of a given particle
at time $t$ is the path describing for any $s<t$ the location of the
particle or whichever ancestor that is alive at time $s$.) We then
create mass uniformly along this
set of particle paths, which then evolves according to the
law $\hat\N_{(\cdot)}$. Loosely speaking, we want to generate a Poisson
random measure with
intensity
\[
\int_0^\infty4\Upsilon_t(dz)\hat
\N_{z_t}(X \in\cdot) \,dt
\]
and add up the resulting measure-valued processes to form $\tilde{X}$.
But we must be careful to represent $\tilde{X}_{D'}$ using only the
portions of $\Upsilon$ corresponding to particles whose historical
paths have never left $D'$.
To formulate this,
for each $t$ and historical path $z$, we define the following map
$X^{t,z}\dvtx\Omega\mapsto\mathcal{M}^{\mathcal{O}_{D-}}$
\[
X^{t,z}(\omega)_{D'}=\cases{X_{D'}(\omega), &\quad
if $\tau_{D'}(z)>t$,
\cr
0, &\quad otherwise.}
\]
Now we generate a Poisson random measure on $\mathcal{M}^{\mathcal
{O}_{D-}}$ with intensity
\[
\mathcal{\lambda}(A)=\int_0^\infty4
\Upsilon_t(dz)\hat\N_{z_t}\bigl(\bigl(X^{t,z}
\bigr)^{-1}(A)\bigr) \,dt,
\]
which is now well defined for any $A\in\tilde{\mathcal{F}}_{\subset
D-}$, since $(X^{t,z})^{-1}(A)$ is $\mathcal{F}_{\subset D-}^{z(t)}$
measurable. (Recall, for any $c\in E$, $\mathcal{F}_{\subset D-}^{c}$
is defined as the $\sigma$-algebra generated by $\{X_{D'}, D'\in
\mathcal
{O}, c\in D'\}$, which is the domain where the measure $\mathbb{N}_c$
is defined.) Adding up the resulting measure-valued processes gives us
$\tilde{X}=(\tilde{X}_{D'})_{D'\in\mathcal{O}_{D-}}$.

More precisely, the $n$-dimensional transition operators for $\tilde
{X}$ with respect to $Q_y^{\nu}$ is given by the formula
\[
Q_y^{I,\nu}(\phi_I):=Q_y^\nu
\bigl(e^{-\langle\tilde{X}_I,\phi
_I\rangle
}\bigr)=Q_y^\nu
\bigl(e^{-\int_0^\infty
4\langle\Upsilon_t,\hat\N_\cdot(1-e^{-\langle X^{t,z}_I,\phi
_I\rangle
})\rangle\,dt} \bigr),
\]
where $I=\{D_1,\ldots,D_k\}\subset\mathcal{O}$, and $y\in D$.

Note that for $y\notin D_I$, $Q_y^{I,\nu}(\phi_I)=Q_y^{I_j,\nu}(\phi
_{I_j})$ for some $j$ such that $y\notin D_j$ (since all paths in the
backbone start from $y$, $X^{t,z}_{D_j}=0$, $\Upsilon_t$ almost all $z$
for all $t$, $Q_y$-almost surely).

The main result of this section is:
%
%
\begin{theorem}\label{theofrag} Assume (\ref{eqregdom}), that $\nu
\in
\mathcal{N}_{0}^{c}$ and that $\gamma_{x,\nu}(y)<\infty$. Then
$\N_y^\nu$-law of $(X_{D'})_{y\in D'\Subset D}$ is the same as the
$Q_y^\nu$-law of $(\tilde{X}_{D'})_{y\in D'\Subset D}$.
\end{theorem}
\begin{pf} Let $y\in U_1\Subset U_2\Subset\cdots\Subset U_i\cdots$ s.t.
$D=\bigcup_{i=1}^{\infty}U_i$. It will suffice to show that $Q^{I,\nu
}_{y}(\phi_I)=N_{y}^{I,\nu}(\phi_I)$ for $I=\{D_1,\ldots,D_k\}$ where
$D_j\subset D_k =U_i$ for a fixed $i$.

Let $D'=U_i$. Define $\Upsilon^{D'}(d\tilde{z})$ as the random measure
on $D'$-valued paths defined by $\Upsilon^{D'}(d\tilde{z})=\lim
_{t\rightarrow\infty}\int\Upsilon_t(dz)1_{\{z^{\tau_{D'}}\in
d\tilde
{z}\}}$ where $z^{\tau_D}$ is the path $z$ stopped at $\tau_{D'}(z)$.
We will write $\Upsilon^{D'}\sim n$ if the support of $\Upsilon^{D'}$
consists of exactly $n$ paths (i.e., exactly $n$ particles of
$\Upsilon$ exit $D'$). Let
\[
Q_y^{I,\nu,n}(\phi_I):=Q_y^\nu
\bigl(e^{-\langle\tilde{X}_I,\phi
_I\rangle}, \Upsilon^{D'}\sim n\bigr).
\]
We will show by induction on $n$ and $\operatorname{card}(I)$ that
%
%
\begin{equation}\label{eqNeqQ}
N_y^{I,\nu,n}(\phi_I)=\gamma_{x,\nu}(y)Q_y^{I,\nu,n}(
\phi_I)
\end{equation}
for all $y\in D_I$. The theorem then follows by summing on $n$.

First observe that, for $\operatorname{card}(I)=1$,
%
%
\begin{equation}
\label{eqQcardI1n}\label{eqQcardI11} \gamma_{x,\nu}(y)Q_y^{I,\nu,n}(
\phi)=\cases{ \gamma_{x,\nu}(y), &\quad$y\notin D', n=1$,
\cr
0, &\quad$y\notin D',n>1$,}
\end{equation}
and for $\operatorname{card}(I)>1$,
%
%
\begin{equation}
\label{eqQcardIk} \gamma_{x,\nu}(y)Q_y^{I,\nu,n}(
\phi_I) = \gamma_{x,\nu}(y)Q_y^{I^j,\nu,n},\qquad
y\notin D_j,j\neq k.
\end{equation}

Comparing these equations to equations (\ref{eqNcardI11}) and (\ref
{eqNcardIk}), we see that equation (\ref{eqNeqQ}) holds on $\partial D'$.

Let $(\xi_t)_{t\geq0}$ be an $\mathcal{L}$-diffusion under the law
$\Pi_y$. $\Pi_y^{4u}$ and $\Pi_y^{4u,\gamma_{x,\nu}}$ will denote the laws
under which $\xi$ is, respectively, an $\mathcal{L}^{4u}$ diffusion, and
$\gamma_{x,\nu}$-transform of an $\mathcal{L}^{4u}$ diffusion. In what
follows $\xi$ will represent the first particle of the branching
backbone which is by construction following a $\gamma_{x,\nu
}$-transform of an $\mathcal{L}^{4u}$-diffusion.
We let $\zeta$ be the lifetime of $\xi$.
Note first that
\begin{eqnarray*}
&&
e^{-\int_0^{t\wedge\tau_{D_I}} 4u(\xi_s) \,ds} e^{-\int
_0^{t\wedge\tau_{D_I}} 4\hat\N_{\xi_s}(1-e^{-\langle X_I,\phi_I
\rangle
})}
\\
&&\qquad=\exp\biggl(-4\int_0^{t\wedge\tau_{D_I}} \bigl[u(
\xi_s)+\N_{\xi
_s}\bigl(e^{-\langle
X_{D'},u\rangle}-e^{-\langle X_I,\phi_I^u \rangle}
\bigr)\bigr] \,ds\biggr)
\\
&&\qquad=\exp-4\int_0^{t\wedge\tau_{D_I}} \bigl[u(
\xi_s)-\N_{\xi
_s}\bigl(1-e^{-\langle
X_{D'},u\rangle}\bigr)+
\N_{\xi_s}\bigl(1-e^{-\langle X_I,\phi_I^u \rangle
}\bigr)\bigr] \,ds
\\
&&\qquad=\mathcal{N}_{t}^{I,D_I}\bigl(\phi_I^u
\bigr).%
\end{eqnarray*}

If $n=1$, then
%
%
\begin{eqnarray}\label{eqQ1}
&&
\gamma_{x,\nu}(y)Q_y^{I,\nu,n}(
\phi_I)
\nonumber
\\
&&\qquad=\gamma_{x,\nu}(y)\Pi_y^{4u,\gamma_{x,\nu}} \bigl(
e^{-\int
_0^{\tau
_{D_I}} 4\hat\N_{\xi_s}(1-e^{-\langle X_I,\phi_I\rangle})\,ds}Q_{\xi
_{\tau
_{D_I}}}^{I,\nu,1}(\phi_I)1_{\tau_{D_I}<\zeta}
\bigr)
\nonumber
\\
&&\qquad=\Pi_y^{4u} \bigl(\gamma_{\nu}(
\xi_{\tau_{D_I}}) e^{-\int_0^{\tau_{D_I}} 4\hat\N_{\xi_s}(1-e^{-\langle
X_I,\phi
_I\rangle
})\,ds}Q_{\xi_{\tau_{D_I}}}^{I,\nu,1}(
\phi_I) \bigr)
\\
&&\qquad=\Pi_y \bigl(\gamma_{\nu}(\xi_{\tau_{D_I}})e^{-\int_0^{\tau_{D_I}}
4u(\xi_s) \,ds}
e^{-\int_0^{\tau_{D_I}} 4\hat\N_{\xi_s}(1-e^{-\langle X_I,\phi
_I\rangle
})\,ds}Q_{\xi_{\tau_{D_I}}}^{I,\nu,1}(\phi_I) \bigr)
\nonumber
\\
&&\qquad=\Pi_y \bigl( \gamma_{\nu}(\xi_{\tau_{D_I}})Q_{\xi_{\tau
_{D_I}}}^{I,\nu,1}(
\phi_I)\mathcal{N}^{I,D_I}_{\tau_{D'}}\bigl(
\phi_I^u\bigr) \bigr).\nonumber
\end{eqnarray}
The first equation is true, because of the strong Markov property of
$\Upsilon_t$ at the first exit time of $D_I$ of its first branch $\xi$.
The second and third equations follow, respectively, from the
definition of $\gamma_{x,\nu}$-transform, and $L^{4u}$ diffusion. If
\mbox{$\operatorname{card}(I)=1$}, equations (\ref{eqQ1}) and (\ref{eqN1}) and
the fact that
$N_y^{I,\nu,1}(\phi_I)=\gamma_{x,\nu}(y)Q_y^{I,\nu,1}(\phi_I)$ on the
boundary of $D_I=D'$ implies that
\[
N_y^{I,\nu,1}(\phi_I)=\gamma_{x,\nu}(y)
Q_y^{I,\nu,1}(\phi_I)
\]
holds for $y\in\bar{D}_I$. If we assume $N_y^{I,\nu,1}(\phi_I)=\gamma
_{x,\nu}(y)Q_y^{I,\nu,1}(\phi_I)$ holds for $y\in\bar{D}_I$ when
$\operatorname{card}(I)=k-1$, then equations (\ref{eqQcardIk}) and
(\ref {eqNcardIk}) and (\ref{eqQ1}) and (\ref{eqN1}) imply that
$N_y^{I,\nu,1}(\phi_I)=\gamma_{x,\nu}(y)Q_y^{I,\nu,1}(\phi_I)$ for
$y\in \bar{D}_I$ when $\operatorname{card}(I)=k$ as well.

If $n\geq2$, then
\[
\gamma_{x,\nu}(y)Q_y^\nu\bigl(e^{-\langle X_I,\phi_I\rangle},
\Upsilon^{D'}\sim n\bigr)=A+B,
\]
where
\begin{eqnarray*}
A&=&\gamma_{x,\nu}(y)Q_y^\nu
\bigl(e^{-\langle X_I,\phi_I\rangle}, \Upsilon^{D'}\sim n\mbox{, the
first branch of
$\Upsilon$ dies inside $D_I$}\bigr),
\\
B&=&\gamma_{x,\nu}(y)Q_y^\nu
\bigl(e^{-\langle X_I,\phi_I\rangle}, \Upsilon^{D'}\sim n\mbox{, the
first branch of
$\Upsilon$ exits $D_I$}\bigr).
\end{eqnarray*}
We compute $A$ as follows:
%
%
\begin{eqnarray}\label{eqQ2A}
A&=&\gamma_{x,\nu}(y)\nonumber\\[-2pt]
&&{}\times\Pi_{y}^{4u,\gamma_{x,\nu}}
\Biggl( 1_{\zeta<\tau_{D_I}}e^{-\int_0^\zeta
4\hat\N_{\xi_s}(1-e^{-\langle X_I,\phi_I\rangle}) \,ds}\nonumber \\[-2pt]
&&\hspace*{53.3pt}{}\times\sum_{m=1}^n
\int\frac{2\gamma_{x,\nu'}(\xi_\zeta)\gamma_{x,\nu-\nu'}(\xi_\zeta
)}{\Gamma_{x,\nu}(\xi_\zeta)}
Q_{\xi_\zeta}^{\nu'}\bigl(e^{-\langle\tilde{X}_I,\phi
_I\rangle},
\Upsilon^{D'}\sim m\bigr)\nonumber\\[-2pt]
&&\hspace*{103.3pt}{}\times Q_{\xi_\zeta}^{\nu-\nu'}
\bigl(e^{-\langle\tilde{X}_I,\phi
_I\rangle}, \Upsilon^{D'}\sim n-m\bigr) K_x
\bigl(\nu;d\nu' \bigr)\Biggr)
\\[-2pt]
&=&\Pi_{y}^{4u} \Biggl(\int_{0}^{\tau_{D_I}}\,dt\,
\Gamma_{x,\nu}(\xi_t)1_{t<\zeta}
e^{-\int_0^t 4\hat\N_{\xi_s}(1-e^{-\langle\tilde{X}_I,\phi
_I\rangle})\,ds} \nonumber\\[-2pt]
&&\hspace*{25pt}{}\times\sum_{m=1}^n\int
K_x\bigl(\nu;d\nu'\bigr)
\frac{2\gamma_{x,\nu'}(\xi_t)\gamma_{x,\nu-\nu'}(\xi_t)}{\Gamma
_{x,\nu}(\xi_t)}\nonumber
\\[-2pt]
&&\hspace*{65pt}{}\times Q_{\xi_t}^{\nu'}\bigl(e^{-\langle\tilde{X}_I,\phi
_I\rangle
},
\Upsilon^{D'}\sim m\bigr)
 Q_{\xi_t}^{\nu-\nu'}\bigl(e^{-\langle\tilde
{X}_I,\phi
_I\rangle},
\Upsilon^{D'}\sim n-m\bigr) \Biggr)
\nonumber
\\[-2pt]
&=&\Pi_y \Biggl( \int_{0}^{\tau_{D_I}}\,dt\,
e^{-\int_0^t 4u(\xi_s) \,ds}e^{-\int_0^t 4\hat\N_{\xi
_s}(1-e^{-\langle
X_I,\phi_I\rangle}) \,ds}\nonumber
\\[-2pt]
&&\hspace*{18pt}{}\times \sum_{m=1}^{n}
\int K_x\bigl(\nu;d\nu'\bigr)
2\gamma_{x,\nu'}(\xi_t)\gamma_{x,\nu-\nu'}(
\xi_t)\nonumber
\\[-2pt]
&&\hspace*{58pt}{}\times
Q_{\xi_t}^{\nu'}
\bigl(e^{-\langle\tilde{X}_I,\phi_I\rangle}, \Upsilon^{D'}\sim m\bigr)
Q_{\xi_t}^{\nu-\nu'}\bigl(e^{-\langle\tilde
{X}_I,\phi
_I\rangle},
\Upsilon^{D'}\sim n-m\bigr) \Biggr)
\nonumber
\\[-2pt]
&=&\sum_{m=1}^{n}\int\Pi_y
\biggl( \int_{0}^{\tau_{D'}} 2 \mathcal
{N}_{t}^{I,D_I}\bigl(\phi_I^u\bigr)
\cdot\gamma_{x,\nu'}(\xi_t)Q_{\xi_t}^{\nu'}
\bigl(e^{-\langle\tilde
{X}_I,\phi
_I}\rangle, \Upsilon^{D'}\sim m\bigr)
\nonumber
\\[-2pt]
&&\hspace*{74.4pt}{}\times\gamma_{x,\nu-\nu'}(\xi_t)Q_{\xi_t}^{\nu
-\nu'}
\bigl(e^{-\langle
\tilde
{X}_I,\phi_I\rangle}, \Upsilon^{D'}\sim n-m\bigr) \,dt \biggr)\nonumber
\\[-2pt]
&&\hspace*{0pt}{}\times
K_x\bigl(\nu;d\nu'\bigr).\nonumber
\end{eqnarray}

The first equation is true, because at the lifetime $\zeta$ of the
first particle $\Upsilon$ branches into\vspace*{-2pt} two new branching
diffusions with joint conditional law $Q_{\xi_\zeta}^{\nu-\nu'}\times
Q_{\xi _\zeta }^{\nu-\nu'}$ given $(\xi_\zeta)_{t\leq\zeta}$, and
$\nu_1$ whose conditional distribution given $(\xi_\zeta)_{t\leq\zeta}$
has density equal to (\ref{eqfragdensity}). The second equation follows
from a well-known fact on $h$-transforms when $h$ is a potential; see,
for example, \citet{SV}, formula 2.2. The third equation follows
from the definition of a killed diffusion.\vadjust{\goodbreak}

We compute $B$ as follows:
%
%
\begin{eqnarray}\label{eqQ2B}\quad
B&=&\gamma_{x,\nu}(y)\Pi_{y}^{4u,\gamma_{x,\nu}}\bigl(
1_{\zeta>\tau
_{D_I}}e^{-\int_{0}^{\tau_{D_I}}4\hat\N_{\xi_s}(1-e^{-\langle
X_I,\phi
_I\rangle})\,ds}Q_{\xi_{\tau_{D_I}}}^{I,\nu,n}(
\phi_I)\bigr)
\nonumber
\\
&=&\Pi_{y}^{4u}\bigl(\gamma_{x,\nu}(
\xi_{\tau_{D_I}})1_{\zeta>\tau
_{D_I}}e^{-\int_{0}^{\tau_{D_I}}4\hat\N_{\xi_s}(1-e^{-\langle
X_I,\phi
_I\rangle})\,ds}Q_{\xi_{\tau_{D_I}}}^{I,\nu,n}(
\phi_I)\bigr)
\nonumber\\[-8pt]\\[-8pt]
&=&\Pi_{y}\bigl(e^{-\int_{0}^{\tau_{D_I}} 4u(\xi_s) \,ds}\gamma_{x,\nu
}(
\xi_{\tau_{D_I}})e^{-\int_{0}^{\tau_{D_I}}4\hat\N_{\xi
_s}(1-e^{-\langle
X_I,\phi_I\rangle})\,ds}Q_{\xi_{\tau_{D_I}}}^{I,\nu,n}(
\phi_I)\bigr)
\nonumber
\\
&=&\Pi_{y} \bigl( \gamma_{\nu}(\xi_{\tau_{D_I}})Q_{\xi_{\tau
_{D_I}}}^{I,\nu,1}(
\phi_I)\mathcal{N}^{I,D_I}_{\tau_{D_I}}\bigl(
\phi_I^u\bigr) \bigr)\nonumber
\end{eqnarray}
by first applying the strong Markov property of $\Upsilon$ at the first
exit of $D_I$ of its first branch, and then again using the definition
of $\gamma_{x,\nu}$ transform and killed diffusion.

We have shown previously that for all $I$,
$N_y^{I,\nu,1}(\phi_I)=\gamma_{x,\nu}(y)Q_y^{I,\nu,1}(\phi_I)$,
$y\in
\bar{D}_I$. Let us assume for all $I$,
$N_y^{I,\nu,m}(\phi_I)=\gamma_{x,\nu}(y)Q_y^{I,\nu,m}(\phi_I)$,
$y\in
\bar{D}_I$, for all $m\leq n-1$. If $\operatorname{card}(I)=1$, comparing
equations (\ref{eqQ2A}) and (\ref{eqQ2B}) with (\ref{eqNn}) and
using the fact that
$N_y^{I,\nu,n}(\phi_I)=\gamma_{x,\nu}(y)Q_y^{I,\nu,n}(\phi_I)$ on the
boundary of $D_I=D'$, we get that
$N_y^{I,\nu,n}(\phi_I)=\gamma_{x,\nu}(y)Q_y^{I,\nu,n}(\phi_I)$
for all
$y\in\bar{D}_I$. If now in addition we assume
$N_y^{I,\nu,n}(\phi_I)=\gamma_{x,\nu}(y)Q_y^{I,\nu,n}(\phi_I)$ holds
for $y\in\bar{D}_I$ when $\operatorname{card}(I)=k-1$, our induction
hypothesis and equations (\ref{eqQcardIk}), (\ref{eqNcardIk}) combined
with equations (\ref{eqQ2A}), (\ref{eqQ2B}) and (\ref{eqNn}) imply
that $N_y^{I,\nu,1}(\phi_I)=\gamma_{x,\nu}(y)Q_y^{I,\nu,1}(\phi_I)$ for
$y\in\bar{D}_I$ when $\operatorname{card}(I)=k$ as well. Hence equation
(\ref{eqNeqQ}) holds for all $I$ and $y\in\bar{D}_I$, and therefore
the proof is complete.
\end{pf}

Above we described the conditional law $\N_y^\nu$. Now we move to an
arbitrary initial measure $\mu\in\mathcal{M}_{D}^{c}$ such that
$H_{x}^{\nu}(\mu)<\infty$, and so need to handle multiple lines of
descent starting from time $0$. In other words, we are going to
describe the distribution of $(X_{D'})_{\mu\in\mathcal
{M}_{D'},D'\Subset D}$, $P_{\mu}^{\nu}$.
Let $\tilde{P}_{\mu}^{\nu}$ be a probability defined on an axillary
probability space $\tilde{\Omega}$ under which a random cluster of
points $(x_i,\nu_i)_{i=1}^{n}$ in $ D\times V^*$ is generated according
to a distribution proportional to
%
%
\begin{equation}\label{eqinitialcluster}
\frac{1}{n!}\bar{K}_{x,n}(\nu;d\nu_1,\ldots,d
\nu_n)\gamma_{x,\nu_1}(x_1)\cdots
\gamma_{x,\nu_n}(x_n)\mu(dx_1)\cdots\mu
(dx_n).
\end{equation}
Note that this makes sense since
\begin{eqnarray*}
&&\sum_{n=1}^{\infty}\frac{1}{n!}
\int_{x_1,\ldots,x_n}\mu(dx_1)\cdots\mu(dx_n)\\
&&\hspace*{3.5pt}\quad{}\times
\int_{\nu_1,\ldots,\nu_n}\bar{K}_{x,n}(\nu;d\nu_1,
\ldots,d\nu_n)
\gamma_{x,\nu_1}(x_1)\cdots
\gamma_{x,\nu_n}(x_n)\\
&&\hspace*{3.5pt}\qquad = H_{x}^{\nu}(\mu)<
\infty.
\end{eqnarray*}
Once the random cluster is generated, corresponding to each point
$(y_i,\nu_i)$ in the cluster, a measure valued process $X^i$ begins to
evolve following a $Q_{y_i}^{\nu_i}$ law independent of everything
else. This\vspace*{1pt} is consistent with our construction of the law $Q_{y_i}^{\nu
_i}$ since almost surely, each $(y_i,\nu_i)$ will satisfy $\gamma_{x,\nu
_i}(y_i)<\infty$ and $\nu_i\in\mathcal{N}_{0}^{c}$. In addition,
independent from all this, another measure valued process $\tilde{X}^0$
evolves following SBM law with spatial motion killed at rate $u$ and
initial measure $\mu$. Let $\tilde{X}=\sum_{i=1}^{n}\tilde{X}^i+X^0$.
%
%
\begin{theorem}\label{theoforest} Assume (\ref{eqregdom}). Let
$H^{\nu
}_x(\mu)<\infty$. Then $P_{\mu}^{\nu}$-law of
\[
(X_{D'})_{\mu\in
\mathcal{M}_{D'},D'\Subset D}
\]
is the same as $\tilde{P}_{\mu}^{\nu
}$-law of $(\tilde{X}_{D'})_{\mu\in\mathcal{M}_{D'},D'\Subset D}$.
\end{theorem}
\begin{pf} Let $P_{\mu}^{I,\nu}$ and $\tilde{P}_{\mu}^{I,\nu}$ denote
the transition operators of $X$ and $\tilde{X}$. That is,
\[
P_{\mu}^{I,\nu}(\phi_I)=P_{\mu}
\bigl(e^{-\langle X_I,\phi_I\rangle}\bigr),
\]
where $I=(D_1,\ldots,D_k)$, $X_I=(X_{D_1},\ldots,X_{D_k})$ and $\phi
_I=(\phi_1,\ldots,\phi_k)$. $\tilde{P}_{\mu}^{I,\nu}(\phi_I)$ is
defined similarly.

Let $\mu$ be compactly supported in $U_1\Subset U_2\Subset\cdots
\Subset
U_i\cdots$ s.t. $D=\bigcup_{i=1}^{\infty}U_i$. It will suffice to show that
$P_{\mu}^{I,\nu}(\phi_I)=\tilde{P}_{\mu}^{I,\nu}(\phi_I)$ for
$I=\{
D_1,\ldots,D_k\}$ where $D_j\in D_k =U_i$ for a fixed $i$.

Recall $P_{\mu}^{I,\nu}=H^{\nu}(\mu)^{-1}P_{\mu}(e^{-\langle
X_I,\phi
_I\rangle}H^{\nu}(X_{D_k}))$. Note
\begin{eqnarray*}
&&
P_\mu\bigl(e^{-\langle
X_{I},\phi_I\rangle}H^{\nu}(X_{D'})
\bigr)
\\
&&\qquad=P_\mu\Biggl(e^{-\langle
X_{I},\phi_I\rangle}\sum_{n=1}^{\infty}
\int\frac{e^{-\langle
X_{D_k},u\rangle}}{n!} \langle X_{D_k},\gamma_{x,\nu_1}\rangle
\cdots\langle X_{D_k},\gamma_{x,\nu_n}\rangle\Biggr)
\\
&&\qquad\quad{}\times\bar{K}_{x,n}(\nu;d\nu_1,\ldots,d
\nu_n)
\\
&&\qquad=\sum_{n=1}^{\infty}\frac{1}{n!}\int
P_\mu\bigl(e^{-\langle
X_{D_k},\phi_I^u\rangle} \langle X_{D_k},
\gamma_{x,\nu_1}\rangle\cdots\langle X_{D_k},\gamma_{x,\nu_n}
\rangle\bigr)
\\
&&\hspace*{39pt}\qquad\quad{}\times\bar{K}_{x,n}(\nu;d\nu_1,\ldots,d
\nu_n).
\end{eqnarray*}
Using the extended moment formula with
\[
l^{I,u}(x)=4\mathbb{N}_{x}\bigl(1-\exp{-\bigl(\langle
X_{D_1},\phi_1\rangle+\cdots+\langle X_{D_k},
\phi_k +u\rangle\bigr)}\bigr),
\]
we rewrite the right-hand side as
\begin{eqnarray*}
&=&\sum_{n=1}^{\infty}\frac{1}{n!}\int
e^{-\langle\mu,l^{I,u}\rangle}\sum_{\pi(n)}\prod_{i=1}^{r}
\biggl\langle\mu,\mathbb{N}_{(\cdot)} \biggl(e^{-\langle
X_{I},\phi_I^u\rangle}
\prod_{j\in C_i}\langle X_{D_k}, \gamma_{x,\nu_j}\rangle\biggr)
\biggr\rangle
\\
&&\hspace*{40pt}\times\bar{K}_{x,n}(\nu;d\nu_1,\ldots,d
\nu_n).
\end{eqnarray*}
By Lemma~\ref{theokerdecompose} we expand the above expression as
\begin{eqnarray*}
&=&e^{-\langle\mu,l^{I,u}\rangle}\sum_{n=1}^{\infty}
\frac
{1}{n!}\sum_{\pi(n)} \int
\bar{K}_{x,r}(\nu;d\nu_1,\ldots,d\nu_r)
\\
&&{}\times\prod_{i=1}^{r}\biggl\langle\mu,\mathbb{N}_{(\cdot)}
\biggl(e^{-\langle
X_{I},\phi_I^u\rangle} \biggl[\int\prod_{j\in C_i}
\langle X_{D_k},\gamma_{x,\nu_{C_i}}\rangle\bar{K}_{x,n_i}(
\tilde{\nu}_i;d\nu_{C_i}) \biggr] \biggr)\biggr\rangle
\\
&=&e^{-\langle\mu,l^{I,u}\rangle
}\sum_{n=1}^{\infty}
\frac{1}{n!}\sum_{r=1}^{n}\int
\bar{K}_{x,r}(\nu;d\nu_1,\ldots,d\nu_r)\sum
_{n_1,\ldots,n_r}\frac
{n!}{n_1!\cdots
n_r!r!}
\\
&&{}\times\prod_{i=1}^{r}
\biggl\langle\mu,\mathbb{N}_{(\cdot)} \biggl(e^{-\langle
X_{D'},u+\phi\rangle} \biggl[\int
\prod_{j\in C_i}\langle X_{D'},
\gamma_{x,\nu_{j}}\rangle\bar{K}_{x,n_i}(\tilde{\nu}_i;d
\nu_{C_i}) \biggr] \biggr)\biggr\rangle
\\
&=&e^{-\langle\mu,l^{I,u} \rangle}\sum_{r=1}^{\infty}
\frac
{1}{r!}\int\bar{K}_{x,r}(\nu;d\nu_1,\ldots,d
\nu_r)
\prod_{i=1}^{r}\sum
_{n_i=1}^{\infty}\frac
{1}{n_i!}\bigl\langle\mu,
\gamma_{x,\tilde{\nu}_i}\mathbb{N}^{\tilde
{\nu
}_i}_{(\cdot)}
\bigl(e^{-\langle
X_{I},\phi_I\rangle} \bigr)\bigr\rangle
\\
&=&e^{-\langle\mu,l^{I,u}\rangle} \sum_{r=1}^{\infty}
\frac
{1}{r!}\int\bar{K}_{x,r}(\nu;d\nu_1,\ldots,d
\nu_r)\prod_{i=1}^{r}\bigl
\langle\mu,\gamma_{x,\tilde{\nu}_i}Q_{(\cdot)}^{\tilde{\nu}_i}
\bigl(e^{-\langle
X_{I},\phi_I\rangle}\bigr)\bigr\rangle.
\end{eqnarray*}

Note that $e^{-\langle\mu,l^{I,u}\rangle}$ is the transition operator
of a SBM whose spatial motion killed at rate $u$, and the rest of the
expression is the transition operator of $\Sigma_{i=1}^{n}
\tilde{X}_i$, where each $\tilde{X}_{i}$ evolves according to
$Q_{y_i}^{\tilde{\nu}_i}$, where the random cluster of points
$(y_i,\nu_i)_{i=1}^{n}$ is selected according to the density
(\ref{eqinitialcluster}). Hence the proof is complete.
\end{pf}

Note that the extended $X$-harmonic functions $H^\nu_x$ were defined by
Dynkin in \citet{Dyn06}, and we have followed this approach throughout.
The results of this section should be viewed as an attempt to clarify
the structure of these extended $X$-harmonic functions, as well as the
structure of the conditioned superprocesses that are obtained from them.

\section*{Acknowledgments}

The authors thank the referees for a very thorough reading, with many
helpful comments.



\printaddresses


\begin{thebibliography}{17}

\bibitem[\protect\citeauthoryear{Bertoin}{2006}]{Bertoin}
\begin{bbook}[mr]
\bauthor{\bsnm{Bertoin},~\bfnm{Jean}\binits{J.}}
(\byear{2006}).
\btitle{Random Fragmentation and Coagulation Processes}.
\bseries{Cambridge Studies in Advanced Mathematics}
\bvolume{102}.
\bpublisher{Cambridge Univ. Press}, \blocation{Cambridge}.
\bid{doi={10.1017/CBO9780511617768}, mr={2253162}}
\bptok{imsref}%
\end{bbook}
\endbibitem

\bibitem[\protect\citeauthoryear{Doob}{1959}]{Doob}
\begin{barticle}[mr]
\bauthor{\bsnm{Doob},~\bfnm{J.~L.}\binits{J.~L.}}
(\byear{1959}).
\btitle{Discrete potential theory and boundaries}.
\bjournal{J. Math. Mech.}
\bvolume{8}
\bpages{433--458; erratum 993}.
\bid{mr={0107098}}
\bptnote{check related}%
\bptok{imsref}%
\end{barticle}
\endbibitem

\bibitem[\protect\citeauthoryear{Dynkin}{2002}]{Dyn02}
\begin{bbook}[mr]
\bauthor{\bsnm{Dynkin},~\bfnm{E.~B.}\binits{E.~B.}}
(\byear{2002}).
\btitle{Diffusions, Superdiffusions and Partial Differential Equations}.
\bseries{American Mathematical Society Colloquium Publications}
\bvolume{50}.
\bpublisher{Amer. Math. Soc.}, \blocation{Providence, RI}.
\bid{mr={1883198}}
\bptok{imsref}%
\end{bbook}
\endbibitem

\bibitem[\protect\citeauthoryear{Dynkin}{2004}]{Dyn04a}
\begin{bbook}[mr]
\bauthor{\bsnm{Dynkin},~\bfnm{E.~B.}\binits{E.~B.}}
(\byear{2004}).
\btitle{Superdiffusions and Positive Solutions of Nonlinear Partial
  Differential Equations}.
\bseries{University Lecture Series}
\bvolume{34}.
\bpublisher{Amer. Math. Soc.}, \blocation{Providence, RI}.
\bid{mr={2089791}}
\bptok{imsref}%
\end{bbook}
\endbibitem

\bibitem[\protect\citeauthoryear{Dynkin}{2006a}]{DynkinIJM06}
\begin{barticle}[mr]
\bauthor{\bsnm{Dynkin},~\bfnm{E.~B.}\binits{E.~B.}}
(\byear{2006}a).
\btitle{A note on {$X$}-harmonic functions}.
\bjournal{Illinois J. Math.}
\bvolume{50}
\bpages{385--394 (electronic)}.
\bid{issn={0019-2082}, mr={2247833}}
\bptok{imsref}%
\end{barticle}
\endbibitem

\bibitem[\protect\citeauthoryear{Dynkin}{2006b}]{Dyn06}
\begin{barticle}[mr]
\bauthor{\bsnm{Dynkin},~\bfnm{E.~B.}\binits{E.~B.}}
(\byear{2006}b).
\btitle{On extreme {$X$}-harmonic functions}.
\bjournal{Math. Res. Lett.}
\bvolume{13}
\bpages{59--69}.
\bid{issn={1073-2780}, mr={2199566}}
\bptok{imsref}%
\end{barticle}
\endbibitem

\bibitem[\protect\citeauthoryear{Etheridge}{1993}]{E93}
\begin{bincollection}[mr]
\bauthor{\bsnm{Etheridge},~\bfnm{Alison~M.}\binits{A.~M.}}
(\byear{1993}).
\btitle{Conditioned superprocesses and a semilinear heat equation}.
In \bbooktitle{Seminar on {S}tochastic {P}rocesses, 1992 ({S}eattle, {WA},
  1992)}.
\bseries{Progress in Probability}
\bvolume{33}
\bpages{89--99}.
\bpublisher{Birkh\"auser}, \blocation{Boston, MA}.
\bid{mr={1278078}}
\bptok{imsref}%
\end{bincollection}
\endbibitem

\bibitem[\protect\citeauthoryear{Etheridge and Williams}{2003}]{EO03}
\begin{barticle}[mr]
\bauthor{\bsnm{Etheridge},~\bfnm{A.~M.}\binits{A.~M.}} \AND
  \bauthor{\bsnm{Williams},~\bfnm{D.~R.~E.}\binits{D.~R.~E.}}
(\byear{2003}).
\btitle{A decomposition of the {$(1+\beta)$}-superprocess conditioned on
  survival}.
\bjournal{Proc. Roy. Soc. Edinburgh Sect. A}
\bvolume{133}
\bpages{829--847}.
\bid{doi={10.1017/S0308210500002699}, issn={0308-2105}, mr={2006204}}
\bptok{imsref}%
\end{barticle}
\endbibitem

\bibitem[\protect\citeauthoryear{Evans}{1993}]{Eva93}
\begin{barticle}[mr]
\bauthor{\bsnm{Evans},~\bfnm{Steven~N.}\binits{S.~N.}}
(\byear{1993}).
\btitle{Two representations of a conditioned superprocess}.
\bjournal{Proc. Roy. Soc. Edinburgh Sect. A}
\bvolume{123}
\bpages{959--971}.
\bid{doi={10.1017/S0308210500029619}, issn={0308-2105}, mr={1249698}}
\bptok{imsref}%
\end{barticle}
\endbibitem

\bibitem[\protect\citeauthoryear{Evans and Perkins}{1990}]{EP90}
\begin{barticle}[mr]
\bauthor{\bsnm{Evans},~\bfnm{Steven~N.}\binits{S.~N.}} \AND
  \bauthor{\bsnm{Perkins},~\bfnm{Edwin}\binits{E.}}
(\byear{1990}).
\btitle{Measure-valued {M}arkov branching processes conditioned on
  nonextinction}.
\bjournal{Israel J. Math.}
\bvolume{71}
\bpages{329--337}.
\bid{doi={10.1007/BF02773751}, issn={0021-2172}, mr={1088825}}
\bptok{imsref}%
\end{barticle}
\endbibitem

\bibitem[\protect\citeauthoryear{Le~Gall}{1999}]{LG}
\begin{bbook}[mr]
\bauthor{\bsnm{Le~Gall},~\bfnm{Jean-Fran{\c{c}}ois}\binits{J.-F.}}
(\byear{1999}).
\btitle{Spatial Branching Processes, Random Snakes and Partial Differential
  Equations}.
\bpublisher{Birkh\"auser}, \blocation{Basel}.
\bid{doi={10.1007/978-3-0348-8683-3}, mr={1714707}}
\bptok{imsref}%
\end{bbook}
\endbibitem

\bibitem[\protect\citeauthoryear{Mselati}{2004}]{M04}
\begin{bbook}[mr]
\bauthor{\bsnm{Mselati},~\bfnm{Beno{\^{\i}}t}\binits{B.}}
(\byear{2004}).
\btitle{Classification and Probabilistic Representation of the Positive
  Solutions of a Semilinear Elliptic Equation}.
\bseries{Mem. Amer. Math. Soc.}
\bvolume{168}
\bpages{xvi+121}.
\bid{issn={0065-9266}, mr={2031708}}
\bptok{imsref}%
\end{bbook}
\endbibitem

\bibitem[\protect\citeauthoryear{Overbeck}{1993}]{O93}
\begin{barticle}[mr]
\bauthor{\bsnm{Overbeck},~\bfnm{L.}\binits{L.}}
(\byear{1993}).
\btitle{Conditioned super-{B}rownian motion}.
\bjournal{Probab. Theory Related Fields}
\bvolume{96}
\bpages{545--570}.
\bid{doi={10.1007/BF01200209}, issn={0178-8051}, mr={1234623}}
\bptok{imsref}%
\end{barticle}
\endbibitem

\bibitem[\protect\citeauthoryear{Overbeck}{1994}]{O94}
\begin{barticle}[mr]
\bauthor{\bsnm{Overbeck},~\bfnm{L.}\binits{L.}}
(\byear{1994}).
\btitle{Pathwise construction of additive {$H$}-transforms of super-{B}rownian
  motion}.
\bjournal{Probab. Theory Related Fields}
\bvolume{100}
\bpages{429--437}.
\bid{doi={10.1007/BF01268988}, issn={0178-8051}, mr={1305781}}
\bptok{imsref}%
\end{barticle}
\endbibitem

\bibitem[\protect\citeauthoryear{Roelly-Coppoletta and Rouault}{1989}]{RR89}
\begin{barticle}[mr]
\bauthor{\bsnm{Roelly-Coppoletta},~\bfnm{Sylvie}\binits{S.}} \AND
  \bauthor{\bsnm{Rouault},~\bfnm{Alain}\binits{A.}}
(\byear{1989}).
\btitle{Processus de {D}awson--{W}atanabe conditionn\'e par le futur lointain}.
\bjournal{C. R. Acad. Sci. Paris S\'er. I Math.}
\bvolume{309}
\bpages{867--872}.
\bid{issn={0764-4442}, mr={1055211}}
\bptok{imsref}%
\end{barticle}
\endbibitem

\bibitem[\protect\citeauthoryear{Salisbury and Verzani}{1999}]{SV}
\begin{barticle}[mr]
\bauthor{\bsnm{Salisbury},~\bfnm{Thomas~S.}\binits{T.~S.}} \AND
  \bauthor{\bsnm{Verzani},~\bfnm{John}\binits{J.}}
(\byear{1999}).
\btitle{On the conditioned exit measures of super {B}rownian motion}.
\bjournal{Probab. Theory Related Fields}
\bvolume{115}
\bpages{237--285}.
\bid{doi={10.1007/s004400050271}, issn={0178-8051}, mr={1720367}}
\bptok{imsref}%
\end{barticle}
\endbibitem

\end{thebibliography}
\end{document}